\documentclass[12pt]{article}
\usepackage[margin=2cm, top=2cm, bottom=2cm]{geometry}
\usepackage[utf8]{inputenc}
\usepackage[T1]{fontenc}
\usepackage{amsmath,amssymb,amsfonts,amsthm,mathtools}
\usepackage{microtype}
\usepackage{graphicx}
\usepackage{subcaption}
\usepackage{booktabs}
\usepackage{multirow}
\usepackage{algorithm}
\usepackage{algorithmic}
\usepackage{xcolor}
\usepackage{url}
\usepackage[hidelinks,hypertexnames=false]{hyperref}
\allowdisplaybreaks
\DeclareGraphicsExtensions{.pdf,.png,.jpg}

\newtheorem{theorem}{Theorem}[section]
\newtheorem{lemma}[theorem]{Lemma}
\newtheorem{proposition}[theorem]{Proposition}

\newtheorem{definition}[theorem]{Definition}
\newtheorem{assumption}[theorem]{Assumption}
\newtheorem{example}[theorem]{Example}

\newtheorem{remark}[theorem]{Remark}

\newcommand{\email}[1]{\protect\href{mailto:#1}{\nolinkurl{#1}}}

\title{Manifold constrained steepest descent for smooth and closed-set optimization\thanks{
\textbf{Funding:} This work is partially supported by a start-up fund at the University of Hong Kong and by the Hong Kong Research Grants Council under ECS project 27301425. A part of the computations was performed using research computing facilities offered by Information Technology Services and the Department of Mathematics at the University of Hong Kong.}}

\author{Kaiwei Yang\thanks{Department of Mathematics, University of Hong Kong,
Hong Kong, China (\email{yangkaiwei8@connect.hku.hk},
\email{lai.lexiao@hku.hk}).}
\and Lexiao Lai\footnotemark[2]}

\hypersetup{
    pdftitle={Manifold constrained steepest descent for smooth and closed-set optimization},
    pdfauthor={Kaiwei Yang and Lexiao Lai},
    pdfkeywords={manifold optimization, linear minimization oracle, preconditioning, projected gradient descent, closed-set optimization, sparse phase retrieval}
}

\begin{document}
\maketitle
\begin{abstract}
We study minimization of smooth functions over feasible sets that have smooth
embedded-manifold structure throughout or only on selected regions, using
linear minimization oracles (LMOs) to determine search directions under
user-chosen norms. Restricting an LMO to a tangent space, however, can require
an iterative inner solve.
We propose \emph{Manifold Constrained Steepest Descent} (MCSD) and a
tangent-projected variant, MCSD-TP, which avoid solving tangent-space LMO subproblems iteratively.
The spectral-norm specialization of MCSD on the Stiefel manifold yields
\emph{SPEL}, which admits an efficient implementation using matrix-sign
computations.
Under standard regularity assumptions, both methods attain an
\(O(\log T/\sqrt T)\) best-iterate stationarity bound. For closed feasible sets, the
hybrid MCSD--PGD combines either smooth method with projected gradient
descent; under the corresponding local smoothness condition, every accumulation
point is Bouligand stationary. Experiments on Stiefel-constrained PCA,
weighted low-rank approximation, and sparse phase retrieval illustrate the
proposed methods.
\end{abstract}

\section{Introduction}\label{sec:intro}
\noindent Let $\mathbb E$ be a Euclidean space. We consider the constrained optimization problem
\begin{equation}\label{eq:closed_problem_intro}
    \min_{x\in C} f(x),
\end{equation}
where \(f:\mathbb E\to\mathbb R\) is smooth and
\(C\subset\mathbb E\) is a nonempty, closed, and generally nonconvex feasible
set with locally smooth regions; smooth manifolds are the special case in
which this local structure is present everywhere. A canonical family is the norm sphere
\(C_{\diamond}(R):=\{x\in\mathbb E:\|x\|_{\diamond}=R\}\), where \(R>0\).
For the Euclidean norm, \(C_{\diamond}(R)\) is a smooth embedded manifold.
For nonsmooth norms, the sphere may instead consist of smooth regions joined
at exceptional points. For example, vector \(\ell_1\)- and
\(\ell_\infty\)-spheres consist of flat faces that meet nonsmoothly, while the
spectral-norm sphere is smooth when its largest singular value is simple and
nonsmooth at singular-value multiplicities
\cite{Lewis1996Spectral,LewisSendov2005}. 
Beyond norm spheres, two representative matrix constraints are the Stiefel
manifold and the determinantal variety. The Stiefel constraint \(X^\top X=I\), used,
for example, in principal component analysis and orthogonality-constrained
learning, defines a smooth embedded manifold
\cite{brockett1991dynamical,edelman1998geometry,huang2018orthogonal}. By
contrast, the determinantal variety
\(\{X:\operatorname{rank}(X)\le r\}\), central to low-rank matrix completion,
is a closed set whose rank-\(r\) stratum is smooth and whose lower-rank points
are singular when \(r<\min\{m,n\}\)
\cite{vandereycken2013low,olikier2026low}. 

Many existing manifold-optimization methods can be viewed as Riemannian
counterparts of Euclidean algorithms: they replace the Euclidean gradient
and, when needed, the Euclidean Hessian with their Riemannian analogues,
compute a tangent-space direction, and apply a retraction. Standard examples
include Riemannian gradient descent, conjugate-gradient, Newton, and
trust-region methods
\cite[Chapters~4 and~6--8]{absil2008optimization}; see also
\cite{edelman1998geometry,absil2007trust}. Riemannian preconditioning follows
the same pattern but selects a metric informed by the Hessian or a structured
approximation, so that the gradient is obtained from a tangent-space
subproblem \cite{mishra2016riemannian}. Such methods have proved effective for
structured low-rank matrix and tensor problems,
with stationary-point and recovery guarantees available in particular
settings \cite{dong2022riemannian,bian2024preconditioned}. Their efficient
implementation, however, relies on structure that makes the inverse metric
and metric-orthogonal tangent projection tractable; otherwise, computing the
direction may require a tangent-space linear system
\cite{mishra2016riemannian,kressner2016preconditioned}. A cheaper ambient
precondition-then-project variant was proposed for structured low-rank
systems, but its convergence was reported only empirically
\cite[Sections~3--5]{kressner2016preconditioned}.

In unconstrained Euclidean optimization, a different way to encode search
geometry is to choose a norm and minimize a linear model over its unit ball;
we refer to this operation as a \emph{linear minimization oracle} (LMO). A prominent
matrix example is \emph{spectral gradient descent} (SpecGD), which obtains its
direction by evaluating the LMO associated with the spectral-norm ball
\cite{carlson2015preconditioned,pethick2025training}. Building on this
viewpoint, the recently proposed Muon optimizer \cite{jordan2024muon} has
been reported to perform strongly in large-language-model
pretraining \cite{liu2025muon,xu2026deepseek}. More generally, these methods
fit within the framework of \emph{Constrained Steepest Descent} (CSD)
\cite{crawshaw2025exploration}, which interprets each update as a
steepest-descent step under a prescribed norm and thereby provides a flexible
mechanism for controlling relative update magnitudes \cite{bernstein2024old}.

Motivated by the success of LMO-based methods, we ask whether computationally
efficient LMO-based algorithms for manifold-constrained problems can also
provide convergence guarantees. Recent work has explored tangent-space LMOs
for Stiefel constraints
\cite{bernstein2025manifolds,kexuefm-11221,buchanan2025mmuonadmm,cesista2025spectralclipping,solonko2026skewon}.
These methods minimize a linear model over the intersection of the tangent
space and a chosen norm ball; see \eqref{eq:manifold_muon}. Earlier proposals
use heuristic or iterative inner solvers
\cite{bernstein2025manifolds,kexuefm-11221,buchanan2025mmuonadmm,cesista2025spectralclipping}.
Li et al.~\cite{li2026intrinsic} instead use Riemannian metrics to induce
unitarily invariant tangent norms and derive closed-form intrinsic LMOs on
several structured matrix manifolds, using a block-product norm on the
Stiefel manifold. For the joint spectral-norm tangent LMO on the Stiefel
manifold, Solonko et
al.~\cite{solonko2026skewon} recently derived an exact closed-form solution
and proposed the efficient Skewon implementation.

Hyperball provides a different route \cite[Algorithm~1]{wen2026fantastic}.
For the Frobenius sphere
\(C_R:=\{x\in\mathbb E:\|x\|_{\mathbb E}=R\}\), its update is
\(x_{t+1}=P_{C_R}(x_t-\alpha_tR\,u_t/\|u_t\|_{\mathbb E})\), where
\(u_t\ne0\) is produced by an arbitrary base optimizer. When Muon (or SpecGD) is the base
optimizer, the direction \(-u_t/\|u_t\|_{\mathbb E}\) is a rescaled ambient
spectral-norm LMO output, so the method has a single outer optimization loop
\cite{pethick2025training,wen2026fantastic}. Projection, however, need not
convert an ambient LMO direction into descent for the constrained objective:
Example~\ref{ex:ambient_lmo} shows that it can eliminate the direction
entirely. Thus, even for smooth constraints, the central difficulty is to
retain flexible and inexpensive LMO directions without losing stationarity
guarantees.

The same difficulty persists for closed sets that are smooth only locally.
Xie et al.~\cite[Sections~3.1--3.4]{xie2026controlled} propose the Spectral
Sphere Optimizer (SSO) and Muon Sphere for solving
\eqref{eq:closed_problem_intro} over the spectral sphere
\(C=\{W:\|W\|_{2\to 2}=R\}\). Their formulation does not address the nonsmooth
points of the sphere or provide a convergence-to-stationarity guarantee. At a
more general level, projected gradient descent (PGD) with backtracking uses
\(x_{t+1}\in P_C(x_t-\alpha_t\nabla f(x_t))\) for an arbitrary nonempty
closed set \(C\). Under continuous differentiability of the objective function, every accumulation point
of the generated sequence, if any, is Bouligand stationary
\cite{OlikierWaldspurger2025PGD,pauwels2024note}.
Structure-aware methods extend the tangent-step-and-return pattern of
Riemannian gradient descent to such sets. In particular,
projected-projected gradient descent (P2GD) projects the negative gradient
onto the tangent cone and then projects the trial point back to the feasible
set
\cite{SchneiderUschmajew2015P2GD,OlikierGallivanAbsil2023Stratified,olikier2026low}.
For the determinantal variety, P2GD--PGD uses P2GD when the nonzero singular
values are separated from zero and PGD near a possible rank drop, while
retaining the Bouligand-stationarity guarantee
\cite[Algorithm~7.2 and Theorem~7.3]{olikier2026low}. These methods retain
Bouligand-stationarity guarantees, but their directions, as formulated, are
defined by the ambient Euclidean gradient rather than by an LMO associated
with a freely chosen norm.

\paragraph{Contributions}
We introduce \emph{Manifold Constrained Steepest Descent} (MCSD), a framework
for manifold optimization that uses LMO-based directions induced by
user-chosen norms while maintaining feasibility. MCSD combines an ambient LMO
with a metric projection back to the manifold, retaining the simplicity of
constrained steepest descent in Euclidean spaces
(Algorithm~\ref{alg:MCSD}). Its spectral-norm specialization on the Stiefel
manifold yields \emph{SPEL}, which can be implemented efficiently using
matrix-sign computations (Table~\ref{tab:mcsd-2x2}). When metric projection is
costly but a retraction is efficient, the MCSD-TP variant instead projects the
LMO direction onto the tangent space before retracting
(Algorithm~\ref{alg:MCSDTP}). Ellipsoidal norms generated by easily
invertible self-adjoint positive-definite operators also let the framework
incorporate (fixed) ambient preconditioning without solving
the corresponding weighted tangent-space system
(Section~\ref{subsec:lmo_prelim}; Remark~\ref{rem:preconditioned_lmo}).
Treating the ambient LMO and return operation as
primitives, both methods have a single outer loop and, under suitable
assumptions, achieve an \(O(\log T/\sqrt T)\) best-iterate stationarity rate
(Theorem~\ref{thm:convergence_deterministic}). For closed feasible sets with
locally smooth structure, we combine either variant with a PGD fallback, using
the smooth update in locally manifold-like regions and PGD elsewhere (Algorithm~\ref{alg:mcsd_pgd}). Under
suitable assumptions, every accumulation point of the generated sequence, if
any, is Bouligand stationary
(Theorem~\ref{thm:mcsd_pgd_bstationary}). Finally, experiments on PCA over the Stiefel manifold, weighted
low-rank approximation over the determinantal variety, and sparse phase
retrieval over sparse vectors illustrate the proposed methods
(Section~\ref{sec:experiments}). Together, these studies emphasize the
practical value of choosing an LMO geometry suited to the problem.

\section{Preliminaries}\label{sec:prelim}

This section introduces the notation and basic facts on manifold optimization, linear
minimization oracles, and variational analysis.

\paragraph{Notation} Throughout, let \(\mathbb E\) be a Euclidean space with inner product
\(\langle\cdot,\cdot\rangle\) and reference norm \(\|\cdot\|_{\mathbb E}\),
and let \(\mathbb N:=\{0,1,2,\ldots\}\). For vectors,
\(\|\cdot\|_p\) denotes the usual \(\ell_p\)-norm. On matrix spaces,
\(\|\cdot\|_F\) denotes the Frobenius norm. For
\(A\in\mathbb R^{m\times n}\) and \(p,q\in[1,\infty]\), the induced
\(p\)-to-\(q\) norm is
\(\|A\|_{p\to q}:=\sup_{x\in\mathbb R^n\setminus\{0\}}
\|Ax\|_q/\|x\|_p\). In particular, \(\|\cdot\|_{2\to 2}\) is the
spectral norm, and \(\|\cdot\|_{\mathbb E}=\|\cdot\|_F\) unless stated
otherwise. For a given norm \(\|\cdot\|\), its dual norm
\(\|\cdot\|_*\) is given by
\[
    \|g\|_* := \sup_{\|s\|\le 1}\langle g,s\rangle
\]
for all \(g\in\mathbb E\); see
\cite[Appendix~A.1.6]{boyd2004convex}.
All derivative bounds and quadratic remainders are measured in
\(\|\cdot\|_{\mathbb E}\). For a linear or multilinear map, the subscript
\({\rm op}\) denotes the induced Euclidean operator norm, and \(A^*\) denotes
the Euclidean adjoint of a linear map \(A\). For a nonempty set
\(\mathcal A\subset\mathbb E\) and \(x\in\mathbb E\), we denote their
Euclidean distance and possibly empty, set-valued nearest-point projection by
\[
    \operatorname{dist}_{\mathbb E}(x,\mathcal A)
    :=\inf_{y\in\mathcal A}\|y-x\|_{\mathbb E},
    \qquad
    P_{\mathcal A}(x)
    :=\operatorname*{argmin}_{y\in\mathcal A}\|y-x\|_{\mathbb E}.
\]
When \(P_{\mathcal A}(x)\) is a singleton, we identify it with its unique
element.

\subsection{Manifold optimization}\label{sec:manifold_opt}

We first specialize the constrained optimization
problem~\eqref{eq:closed_problem_intro} to \(C=\mathcal M\), where
\(\mathcal M\subset\mathbb E\) is a \(C^2\) embedded manifold without boundary,
equipped with the Riemannian metric induced by the ambient inner product
\cite[Definitions~3.10 and~3.55]{boumal2023introduction}. We write
\(T_x\mathcal M\) for the tangent space at \(x\). If the objective \(f\) is
\(C^1\), then its Riemannian gradient is the tangent projection of the ambient
gradient,
\[
    \nabla_{\mathcal M} f(x)
    := P_{T_x\mathcal M}(\nabla f(x)).
\]
See \cite[Definition~3.58 and Proposition~3.61]{boumal2023introduction}.

To describe the return maps used below, we first recall the tangent bundle and
fix our derivative notation. The tangent bundle \cite[Definition~3.42]{boumal2023introduction} of \(\mathcal M\) is the
disjoint union
\[
    T\mathcal M
    :=\{(x,\xi):x\in\mathcal M,\ \xi\in T_x\mathcal M\}.
\]
When the base point is clear, we identify a tangent vector
\(\xi\in T_x\mathcal M\) with \((x,\xi)\in T\mathcal M\).

For a \(C^1\) map \(F:\mathcal N\to\mathcal P\) between \(C^1\) manifolds,
its differential at \(z\) is the linear map
\(DF(z):T_z\mathcal N\to T_{F(z)}\mathcal P\) defined by
\[
    DF(z)[v]
    :=\left.\frac{d}{dt}\right|_{t=0}F(\gamma(t)),
\]
where
\(\gamma:(-\varepsilon,\varepsilon)\to\mathcal N\) is a \(C^1\) curve with
\(\gamma(0)=z\) and \(\gamma'(0)=v\); see
\cite[Definition~3.34]{boumal2023introduction}.

With this notation in place, a \emph{retraction} on \(\mathcal M\) is a
\(C^1\) map
\[
    \operatorname{Retr}:T\mathcal M\to\mathcal M,
    \qquad
    (x,\xi)\mapsto\operatorname{Retr}_x(\xi),
\]
such that, for every \((x,\xi)\in T\mathcal M\), the curve
\(c(t):=\operatorname{Retr}_x(t\xi)\) satisfies \(c(0)=x\) and
\(c'(0)=\xi\); see
\cite[Definition~3.47]{boumal2023introduction}. A retraction therefore
provides a first-order accurate return from \(T_x\mathcal M\) to
\(\mathcal M\). A \emph{local retraction} is a \(C^1\) map defined on an
open neighborhood of the zero vectors in \(T\mathcal M\) and satisfying the
same curve conditions
\cite[discussion following Definition~3.47]{boumal2023introduction}.
Classical Riemannian gradient descent (RGD) uses
\(x_{t+1}=\operatorname{Retr}_{x_t}
(-\alpha_t\nabla_{\mathcal M}f(x_t))\) with \(\alpha_t>0\).

A canonical \(C^\infty\) embedded manifold used in this paper is the Stiefel
manifold \cite[Section~7.3]{boumal2023introduction}. For \(1\le p\le n\), it
is given by
\[
    \mathrm{St}(n,p)
    := \left\{ X \in \mathbb R^{n\times p} : X^\top X = I_p \right\},
\]
embedded in $\mathbb R^{n\times p}$. It reduces to the unit sphere when $p = 1$. For any $X\in\mathrm{St}(n,p)$, its tangent space is
\[
    T_X \mathrm{St}(n,p)
    = \left\{
        Z \in \mathbb R^{n\times p}
        :
        X^\top Z + Z^\top X = 0
      \right\}.
\]
The orthogonal projection of $G\in\mathbb R^{n\times p}$ onto $T_X\mathrm{St}(n,p)$ is
\[
    P_{T_X\mathrm{St}(n,p)}(G)
    \;=\; G - X\,\mathrm{sym}(X^\top G),
\]
where $\mathrm{sym}(A):=\tfrac12(A+A^\top)$. For a nonzero
$G\in\mathbb R^{n\times p}$ of rank $r$, let
$G=U_r\Sigma_rV_r^\top$ be a compact singular value decomposition. We define
the \emph{generalized matrix sign}, or partial polar factor, by
\begin{equation}\label{eq:msign}
    \mathrm{msign}(G)
    :=U_rV_r^\top
    =G\bigl((G^\top G)^\dagger\bigr)^{1/2},
    \qquad
    \mathrm{msign}(0):=0,
\end{equation}
where ${}^\dagger$ denotes the Moore--Penrose pseudoinverse. This definition
is independent of the chosen compact SVD. If
$Y\in\mathbb R^{n\times p}$ has full column rank, then
\[
    P_{\mathrm{St}(n,p)}(Y)
    =\mathrm{msign}(Y)
    =Y(Y^\top Y)^{-1/2}.
\]
Thus, at full column rank, $\mathrm{msign}(Y)$ is the unique Euclidean
projection of $Y$ onto $\mathrm{St}(n,p)$. At deficient rank,
\eqref{eq:msign} remains well defined but is generally not an element of the
Stiefel manifold.

\subsection{Linear minimization oracles}
\label{subsec:lmo_prelim}

The \emph{linear minimization oracle} (LMO) associated with a given norm \(\|\cdot\|\)
is the set-valued map
\[
    \mathrm{LMO}_{\|\cdot\|}(g)
    :=
    \operatorname*{argmin}_{\|s\|\le 1}\langle g,s\rangle,
    \qquad g\in\mathbb E.
\]
This is the usual linear-optimization oracle specialized to the norm unit
ball; see \cite[Sections~1.4 and~2.1]{braun2025conditional}.
Its optimal value is \(-\|g\|_*\). Hence every
\(s\in\mathrm{LMO}_{\|\cdot\|}(g)\) satisfies
\(\langle g,s\rangle=-\|g\|_*\). Throughout, an iteration-indexed LMO
output is denoted by \(s_t\); any subsequent transformation of that output is
denoted separately.
Such oracles produce steepest-descent directions for the geometry induced by
the prescribed norm; see
\cite{jaggi2013revisiting,crawshaw2025exploration,pethick2025training}.

An important special case is a fixed \emph{ellipsoidal}, or
\emph{preconditioned}, norm. Let
\(\mathcal B:\mathbb E\to\mathbb E\) be self-adjoint and positive definite,
and define
\begin{equation}
\label{eq:preconditioned_lmo}
    \|u\|_{\mathcal B}:=\langle u,\mathcal Bu\rangle^{1/2},
    \qquad
    \|g\|_{\mathcal B,*}:=\langle g,\mathcal B^{-1}g\rangle^{1/2}.
\end{equation}
For \(g\ne0\), the associated LMO is single-valued and satisfies
\[
    \mathrm{LMO}_{\|\cdot\|_{\mathcal B}}(g)
    =\left\{-\frac{\mathcal B^{-1}g}{\|g\|_{\mathcal B,*}}\right\};
\]
at \(g=0\), we select the output zero. Thus the oracle applies
\(\mathcal B^{-1}\) in the ambient space and then normalizes the result;
when \(\mathcal B\) models curvature, \(\mathcal B^{-1}\) acts as an ambient
preconditioner.

When $\mathbb E$ is a matrix space, different matrix norms lead to different closed-form LMOs; see, e.g., \cite[Tables 2--4]{pethick2025training}.
For instance, if $\mathbb E=\mathbb R^{n\times p}$ and $\|\cdot\| = \|\cdot\|_{2\to 2}$, then one may take
\[
    -\mathrm{msign}(G)\in\mathrm{LMO}_{\|\cdot\|_{2\to 2}}(G)
\]
for every $G\in \mathbb R^{n\times p}$. Indeed, for $G\ne0$,
$\|\mathrm{msign}(G)\|_{2\to 2}=1$ and
\[
    \langle G,-\mathrm{msign}(G)\rangle
    =-\operatorname{tr}(\Sigma_r)
    =-\|G\|_*,
\]
while the claim is immediate from $\mathrm{msign}(0)=0$ when $G=0$.
This choice yields the update direction used in spectral gradient descent \cite{carlson2015preconditioned}
and, more recently, in the Muon optimizer \cite{jordan2024muon}.
These algorithms admit cost-effective implementations, thanks to the GPU-friendly iterative schemes for approximating $\mathrm{msign}(\cdot)$, including the
classical Newton--Schulz iteration \cite{schulz1933iterative} and recent accelerated variants; see, e.g.,
\cite{grishina2025accelerating,amsel2025polar}. 

\subsection{Variational analysis}\label{sec:closed_set_prelim}

We also consider \eqref{eq:closed_problem_intro} with \(C\subset\mathbb E\)
nonempty and closed. Following standard variational-analysis notation
\cite[Definition~6.1 and Proposition~6.5]{RockafellarWets1998}, for \(x\in C\) the
Bouligand tangent cone is
\[
    T_C(x):=
    \left\{
        v\in\mathbb E:
        \exists \tau_j\downarrow0,\ \exists v_j\to v
        \text{ such that } x+\tau_jv_j\in C
    \right\}.
\]
The regular normal cone is
\[
    \widehat N_C(x):=
    \left\{
        w\in\mathbb E:
        \langle w,v\rangle\le0
        \text{ for all } v\in T_C(x)
    \right\}.
\]
A point \(x\in C\) is \emph{Bouligand stationary} \cite[Definition~6.1.1]{cui2021modern}, or \emph{B-stationary},
for \(\min_{x\in C}f(x)\) if $-\nabla f(x)\in\widehat N_C(x)$.

\section{Algorithms}\label{sec:MCSD}

This section develops three methods for
problem~\eqref{eq:closed_problem_intro}.
Subsection~\ref{subsec:mcsd_variants} treats the smooth manifold case
\(C=\mathcal M\), introducing manifold constrained steepest descent with a
metric-projection return and its tangent-projected retraction variant.
Subsection~\ref{sec:mcsd_pgd} combines either method with projected gradient
descent for a general closed set \(C\).

\subsection{Manifold constrained steepest descent methods}
\label{subsec:mcsd_variants}

Consider \eqref{eq:closed_problem_intro} with \(C=\mathcal M\), where
\(\mathcal M\subset\mathbb E\) is a \(C^2\) embedded manifold, and suppose that the
metric projection \(P_{\mathcal M}\) can be computed efficiently. Starting
from \(x_t\in\mathcal M\) with \(\alpha_t>0\), a projection-return step based
on a search direction \(s_t\in\mathbb E\) has the form
\begin{equation}\label{eq:proj_update}
    x_{t+1}
    = P_{\mathcal M}\bigl(x_t + \alpha_t s_t\bigr),
    \qquad \forall t \in \mathbb N.
\end{equation}
A seemingly natural choice is
\(s_t\in\mathrm{LMO}_{\|\cdot\|}(\nabla f(x_t))\), but even an arbitrarily
small step of this form need not decrease \(f\).

\begin{example}[An ambient-gradient LMO need not descend]\label{ex:ambient_lmo}
Let \(\mathcal M=\mathbb S^1\subset\mathbb R^2\), \(x=e_1\),
\(f(z)=2z_1+z_2\), and \(\|\cdot\|=\|\cdot\|_1\). Since
\(\nabla f(x)=(2,1)\), the unique LMO output is \(s=-e_1\). Yet
\(\nabla_{\mathcal M}f(x)=(I-xx^\top)\nabla f(x)=e_2\ne0\), while for every
\(0<\alpha<1\), update~\eqref{eq:proj_update} gives
\(P_{\mathbb S^1}(x+\alpha s)=P_{\mathbb S^1}((1-\alpha)e_1)=x\).
Hence \(f(P_{\mathbb S^1}(x+\alpha s))=f(x)\): the nonstationary point \(x\)
is fixed by every sufficiently small step.
\end{example}

The obstruction is geometric: the
ambient-gradient LMO selects a normal direction, which projection removes.
The appropriate oracle input is revealed by linearizing the projected
objective. Locally around \(x_t\in\mathcal M\), the nearest-point projection
is single-valued and \(C^1\), and
\(\nabla(f\circ P_{\mathcal M})(x_t)=\nabla_{\mathcal M}f(x_t)\)
\cite[Lemma~4]{absil2012projection}. If
\(\nabla(f\circ P_{\mathcal M})\) is locally \(L_P\)-Lipschitz, the Euclidean
descent lemma \cite[Lemma~1.2.3]{nesterov2003introductory} gives
\begin{equation}\label{eq:descent}
    \begin{aligned}
    f(x_{t+1}) &= f\circ P_{\mathcal M}(x_t + \alpha_t s_t)\\
    &\le f(x_t) + \alpha_t \langle \nabla_{\mathcal M}f(x_t),s_t\rangle
    + \frac{L_P\alpha_t^2}{2}\|s_t\|^2_\mathbb{E}.
    \end{aligned}
\end{equation}
If \(\nabla_{\mathcal M}f(x_t)\ne0\), choosing
\(s_t\in\mathrm{LMO}_{\|\cdot\|}
(\nabla_{\mathcal M}f(x_t))\) makes the linear term equal to
\(-\|\nabla_{\mathcal M}f(x_t)\|_*\), which dominates the quadratic
remainder for sufficiently small \(\alpha_t\). This yields the
\emph{manifold constrained steepest descent} (MCSD) method in
Algorithm~\ref{alg:MCSD}.

\begin{algorithm}[!ht]
\caption{Manifold constrained steepest descent (MCSD)}
\label{alg:MCSD}
\begin{algorithmic}[1]
\REQUIRE Initial point \(x_0\in\mathcal M\), a norm \(\|\cdot\|\) and its
LMO, and stepsizes \((\alpha_t)_{t\in\mathbb N}\).
\FOR{\(t=0,1,2,\ldots\)}
    \IF{\(\|\nabla_{\mathcal M}f(x_t)\|_*=0\)}
        \STATE \textbf{stop} and return \(x_t\).
    \ENDIF
    \STATE Choose
    \(s_t\in\mathrm{LMO}_{\|\cdot\|}
    \bigl(\nabla_{\mathcal M}f(x_t)\bigr)\).
    \STATE \(x_{t+1}=P_{\mathcal M}(x_t+\alpha_ts_t)\).
\ENDFOR
\end{algorithmic}
\end{algorithm}

MCSD extends the constrained steepest descent framework
\cite{crawshaw2025exploration} (see also \cite{pethick2025training}) from
Euclidean spaces to smooth embedded manifolds.
Table~\ref{tab:mcsd-2x2} lists representative instances. Its spectral-norm
specialization on the Stiefel manifold is SPEL, which, to the best of our
knowledge, has not appeared in the existing literature. This specialization is
of interest because its per-iteration cost is only slightly higher than that of
RGD---one additional \(\mathrm{msign}\) operation---and much lower than that
of other Manifold Muon approaches
\cite{kexuefm-11215,bernstein2025manifolds,buchanan2025mmuonadmm}.

\begin{table*}[!ht]
\centering
\caption{Representative one-step MCSD updates for choices of the underlying
manifold \(\mathcal M\) and LMO norm. Here, \(\|\cdot\|_{2\to 2}\) denotes the matrix
spectral norm, and
SPEL abbreviates \emph{SP}ectral gradient descent on the Stief\emph{EL}
manifold. NormGD and NormRGD denote normalized gradient descent and normalized
Riemannian gradient descent, respectively, while SpecGD denotes spectral
gradient descent. The Stiefel formulas are understood locally where the trial
matrix has full column rank.
}
\label{tab:mcsd-2x2}
\renewcommand{\arraystretch}{1.25}
\setlength{\tabcolsep}{5pt}
\resizebox{\textwidth}{!}{
\begin{tabular}{@{}lcc@{}}
\toprule
& \(\|\cdot\|=\|\cdot\|_{\mathbb E}\)
& \(\|\cdot\|=\|\cdot\|_{2\to 2}\) \\
\midrule
\(\mathcal M=\mathbb R^{n\times p}\)
&
\textup{NormGD:} 
\(\displaystyle
X_{t+1}
=
X_t-\alpha_t
\frac{\nabla f(X_t)}
{\|\nabla f(X_t)\|_{\mathbb E}}
\)
&
\textup{SpecGD:}
\(\displaystyle
X_{t+1}
=
X_t-\alpha_t\,
\mathrm{msign}\!\bigl(\nabla f(X_t)\bigr)
\)
\\
\addlinespace[0.5em]
\(\mathcal M=\mathrm{St}(n,p)\)
&
\textup{NormRGD:}
\(\displaystyle
X_{t+1}
=
\mathrm{msign}\!\left(
X_t-\alpha_t
\frac{\nabla_{\mathcal M} f(X_t)}
{\|\nabla_{\mathcal M} f(X_t)\|_{\mathbb E}}
\right)
\)
&
\textup{SPEL:}
\(\displaystyle
X_{t+1}
=
\mathrm{msign}\!\left(
X_t-\alpha_t\,
\mathrm{msign}\!\bigl(\nabla_{\mathcal M}f(X_t)\bigr)
\right)
\)
\\
\bottomrule
\end{tabular}
}
\end{table*}

For a given manifold constraint $\mathcal M$ and norm $\|\cdot\|$,
an efficient implementation of Algorithm~\ref{alg:MCSD} requires efficient
computations of both an element of $\mathrm{LMO}_{\|\cdot\|}$ and an
element of the metric projection $P_{\mathcal M}$. For common norms, the LMO
is often available in closed form, reducing to coordinate selection or leading
singular/eigenvector computations; see, e.g.,
\cite{jaggi2013revisiting,pethick2025training}. Metric projections are also
tractable, at least locally, on many matrix manifolds: examples include the
Stiefel manifold via $\mathrm{msign}$ (see~\eqref{eq:msign}), the fixed-rank
manifold via truncated SVD \cite{vandereycken2013low}, and suitable symmetric
spectral manifolds via eigenvalue projections
\cite{absil2012projection,lewis2008alternating}.

Outside these favorable cases, however, the metric projection may be
unavailable in closed form or expensive to evaluate. A retraction instead acts only on tangent steps and needs to agree
with them only to first order
\cite[Section~5.12 and Definition~3.47]{boumal2023introduction}; it can also be substantially easier to compute. On the open
positive-definite cone, the quadratic retraction in
\cite[Eq.~(4.10)]{jeuris2012survey} is defined for every tangent vector,
whereas an indefinite or singular trial matrix has no nearest point in the
open cone. Likewise, on the fixed-multilinear-rank tensor manifold, truncated
HOSVD is a computable local retraction that avoids computing the exact metric
projection
\cite[Proposition~3 and Section~3.3]{kressner2014lowranktensor}. To exploit
such retractions, we
project the ambient LMO output $s_t$ onto the tangent space,
$d_t=P_{T_{x_t}\mathcal M}s_t$, and set
$x_{t+1}=\operatorname{Retr}_{x_t}(\alpha_t d_t)$. The resulting algorithm is summarized in Algorithm~\ref{alg:MCSDTP}.

\begin{algorithm}[!ht]
\caption{Tangent-projected MCSD (MCSD-TP)}
\label{alg:MCSDTP}
\begin{algorithmic}[1]
\REQUIRE Initial point \(x_0\in\mathcal M\), a norm \(\|\cdot\|\) and its
LMO, step sizes \((\alpha_t)_{t\in\mathbb N}\), and a \(C^2\) retraction
\(\operatorname{Retr}\) on \(\mathcal M\).
\FOR{\(t=0,1,2,\ldots\)}
    \IF{\(\|\nabla_{\mathcal M}f(x_t)\|_*=0\)}
        \STATE \textbf{stop} and return \(x_t\).
    \ENDIF
    \STATE Choose
    \(s_t\in\mathrm{LMO}_{\|\cdot\|}
    (\nabla_{\mathcal M}f(x_t))\).
    \STATE Set
    \(d_t=P_{T_{x_t}\mathcal M}s_t\).
    \STATE
    \(x_{t+1}=\operatorname{Retr}_{x_t}(\alpha_t d_t)\).
\ENDFOR
\end{algorithmic}
\end{algorithm}

\begin{remark}[Choosing between MCSD and MCSD-TP]
Unlike MCSD, MCSD-TP does not require the metric projection
$P_{\mathcal M}$: as its return map, it requires only a \(C^2\) retraction
after the LMO output has been projected onto the tangent space. It is therefore
applicable when a suitable retraction and tangent projection are available,
even if $P_{\mathcal M}$ is unavailable or impractical. When
$P_{\mathcal M}$ can be computed efficiently, however, we recommend the
simpler MCSD update; MCSD-TP is intended primarily for settings in which a
retraction provides a more practical return map.
\end{remark}

\begin{remark}[Quotient-norm interpretation of MCSD-TP]
\label{rem:quotient_norm}
Fix \(x\in\mathcal M\), define
\[
    \|d\|_{x,\mathrm q}
    :=
    \inf\{\|s\|:P_{T_x\mathcal M}s=d\},
    \qquad d\in T_x\mathcal M.
\]
This is a norm on \(T_x\mathcal M\), whose unit ball is the projected
ambient unit ball \(P_{T_x\mathcal M}\{s:\|s\|\leq1\}\). Since
\(g:=\nabla_{\mathcal M}f(x)\in T_x\mathcal M\),
\[
    \min_{\|d\|_{x,\mathrm q}\leq1}\langle g,d\rangle
    =
    \min_{\|s\|\leq1}\langle g,P_{T_x\mathcal M}s\rangle
    =
    \min_{\|s\|\leq1}\langle g,s\rangle
    =
    -\|g\|_*.
\]
Consequently, if \(s\in\mathrm{LMO}_{\|\cdot\|}(g)\), then \(P_{T_x\mathcal M}s\)
is a steepest-descent direction under the quotient norm
\(\|\cdot\|_{x,\mathrm q}\). Thus MCSD-TP realizes a tangent-space
steepest-descent step through an ambient LMO followed by an orthogonal
projection. 
\end{remark}

\begin{remark}[Relation to Riemannian preconditioning]
\label{rem:preconditioned_lmo}
Specializing the preceding quotient-norm interpretation to the ellipsoidal
norm \(\|\cdot\|_{\mathcal B}\) in \eqref{eq:preconditioned_lmo} clarifies
the relation between MCSD-TP and Riemannian preconditioning. Let \(g_t:=\nabla_{\mathcal M}f(x_t)\ne0\). Then the MCSD-TP direction is
\[
    d_t=P_{T_{x_t}\mathcal M}s_t
    =-\frac{P_{T_{x_t}\mathcal M}\mathcal B^{-1}g_t}{\|g_t\|_{\mathcal B,*}},
    \qquad
    \langle g_t,d_t\rangle=-\|g_t\|_{\mathcal B,*},
\]
where the second identity follows from \(g_t=P_{T_{x_t}\mathcal M}g_t\). By contrast, under
the Riemannian metric
\(\langle u,v\rangle_{\mathcal B}:=\langle u,\mathcal Bv\rangle\), the
Riemannian gradient is the unique \(v_t\in T_{x_t}\mathcal M\) satisfying
\[
    \left.(P_{T_{x_t}\mathcal M}\mathcal BP_{T_{x_t}\mathcal M})\right|_{T_{x_t}\mathcal M}v_t=g_t.
\]
Thus \(d_t\) is generally not the normalized negative
\(\mathcal B\)-metric Riemannian gradient \cite{mishra2016riemannian,bian2024preconditioned}: the latter requires a weighted
tangent-space solve, whereas MCSD-TP applies \(\mathcal B^{-1}\) in the
ambient space and then uses the ordinary Euclidean tangent projection.  
\end{remark}

\subsection{A hybrid algorithm for closed-set optimization}
\label{sec:mcsd_pgd}

We now return to the general closed-set problem
\eqref{eq:closed_problem_intro}, where \(C\subset\mathbb E\) is nonempty and
closed and \(f:\mathbb E\to\mathbb R\) is \(C^1\). Because \(C\) may contain
nonsmooth points, manifold optimization methods are not applicable everywhere. A natural
remedy is to use a manifold optimization method on smooth pieces and the projected gradient descent (PGD) elsewhere;
PGD provides stationarity guarantees under minimal assumptions
\cite{OlikierWaldspurger2025PGD,pauwels2024note}. P2GD--PGD follows this
strategy on the determinantal variety, exploiting its stratified structure
for cheaper updates on smooth pieces while invoking PGD near nonsmooth points
\cite{olikier2026low}.

To extend the flexible norm geometry of MCSD to closed-set optimization, we
adopt the same strategy. We first represent the two smooth updates in
Algorithms~\ref{alg:MCSD} and~\ref{alg:MCSDTP} through a common map.

\begin{definition}[Ambient return map]
\label{def:ambient_return_map}
Let \(\mathcal M\subset\mathbb E\) be a \(C^2\) embedded manifold and let
\(x\in\mathcal M\). An \emph{ambient return map of \(\mathcal M\) at \(x\)}
is a \(C^2\) map \(\mathcal R_x:\mathcal D_x\to\mathcal M\), where
\(\mathcal D_x\subset\mathbb E\) is an open neighborhood of \(0\), such that
\[
    \mathcal R_x(0)=x,
    \qquad
    D\mathcal R_x(0)=P_{T_x\mathcal M}.
\]
Here and below, derivatives of manifold-valued maps are taken after viewing
the maps as \(\mathbb E\)-valued.

For \(A\subset\mathcal M\), an \emph{ambient return map on \(A\)} is a
\(C^2\) map \(\mathcal R:\mathcal D\to\mathcal M\), where
\(\mathcal D\subset\mathcal M\times\mathbb E\) is open and contains
\(A\times\{0\}\), such that, for every \(x\in A\), the fiber
\(\mathcal R_x:\mathcal D_x\to\mathcal M\), defined by
\(\mathcal D_x:=\{u\in\mathbb E:(x,u)\in\mathcal D\}\) and
\(\mathcal R_x(u):=\mathcal R(x,u)\), is an ambient return map of
\(\mathcal M\) at \(x\).
\end{definition}
The projection and tangent-projected retraction updates correspond,
respectively, to
\begin{equation}\label{eq:return_maps}
     \mathcal R_x^{\rm P}(u):=P_{\mathcal M}(x+u),
    \qquad
    \mathcal R_x^{\rm R}(u)
    :=\operatorname{Retr}_x(P_{T_x\mathcal M}u),   
\end{equation}
on their respective local domains. The derivative identity follows from
\cite[Lemma~4]{absil2012projection} for \(\mathcal R^{\rm P}\) and from the
defining property of a retraction
for \(\mathcal R^{\rm R}\)
\cite[Definition~3.47]{boumal2023introduction}. To apply these smooth updates to
\(C\), we isolate a relatively closed region on which \(C\) is locally
smooth.

\begin{definition}[Safe region with return maps]
\label{def:safe_region}
Let \(C\subset\mathbb E\) be nonempty and closed. A \emph{safe region with
return maps} for \(C\) is a pair \((\Omega,\mathcal R)\), where
\(\mathcal R:=(\mathcal R_x:\mathcal D_x\to C)_{x\in\Omega}\) and each
\(\mathcal D_x\subset\mathbb E\) is an open neighborhood of \(0\), such that
\begin{enumerate}
    \item \(\Omega\subset C\) is relatively closed in \(C\);
    \item for every \(\bar x\in\Omega\), there exist an open neighborhood
    \(U_{\bar x}\subset\mathbb E\) of \(\bar x\), a \(C^3\) embedded
    manifold \(\mathcal M_{\bar x}\subset\mathbb E\) such that
    \[
        C\cap U_{\bar x}
        =\mathcal M_{\bar x}\cap U_{\bar x}.
    \]
    Moreover, there is an ambient return map of
    \(\mathcal M_{\bar x}\) on \(\Omega\cap U_{\bar x}\) whose fiber at
    every \(x\in\Omega\cap U_{\bar x}\) is a restriction of
    \(\mathcal R_x\).
\end{enumerate}
We call \(\Omega\) the safe region and \(\mathcal R\) its family of
return maps.
\end{definition}

Intuitively, \(\Omega\) selects a relatively closed portion of \(C\) that
stays within its locally smooth part and away from junctions where additional
branches meet. Although \(C\) need not be a manifold globally, its local smooth
structure at each \(x\in\Omega\) is unambiguous: any two local manifold
representatives coincide near \(x\), so \(T_C(x)\) is a well-defined linear
space whose dimension is locally constant, and hence constant on each
connected component of \(\Omega\).
Either canonical map in \eqref{eq:return_maps}, restricted to remain in the
corresponding local manifold, supplies such a local family. The required
joint \(C^2\) dependence is only local; no compatibility is imposed between
extensions associated with distinct neighborhoods. Thus the definition depends
only on the local geometry of \(C\), not on \(f\), and, unlike the related
framework in \cite[Assumption~2.2]{OlikierGallivanAbsil2023Stratified}, does
not require an ordered stratification of all of \(C\).

Every closed semialgebraic set, that is, every closed set defined by finite
Boolean combinations of polynomial equalities and inequalities
\cite[Section~2.1]{bochnak2013real}, admits a finite \(C^3\) semialgebraic
stratification \cite[Proposition~9.1.8]{bochnak2013real}. In particular, a
subset of \(C\) that is relatively closed in \(C\) and consists of points
where \(C\) itself agrees locally with a smooth stratum can be equipped with
return maps to form a safe region with return maps. For many practical
constraints, these objects can be constructed explicitly; see
Examples~\ref{ex:sparse_safe_region}--
\ref{ex:spectral_sphere_safe_region} in
Section~\ref{sec:examples_safe_regions}.

Given a safe region with return maps \((\Omega,\mathcal R)\),
Algorithm~\ref{alg:mcsd_pgd} combines a smooth-region branch with a PGD
fallback. At \(x_t\in\Omega\), it uses the intrinsic tangent space
\(T_C(x_t)\) and applies \(\mathcal R_{x_t}\). Outside \(\Omega\), it applies
PGD to the full set \(C\). Backtracking selects \(\alpha_t\) using
\eqref{eq:mcsd_pgd_smooth_branch} or \eqref{eq:pgd_fallback}, respectively.

\begin{algorithm}[H]
\caption{MCSD--PGD}
\label{alg:mcsd_pgd}
\begin{algorithmic}[1]
\REQUIRE \(x_0\in C\); a norm \(\|\cdot\|\) and its LMO; safe region with
return maps \((\Omega,\mathcal R)\); \(\bar\alpha_M,\bar\alpha_P>0\);
\(c_M,c_P,\beta_M,\beta_P\in(0,1)\).
\FOR{\(t=0,1,2,\ldots\)}
    \IF[Stationarity test]{\(-\nabla f(x_t)\in\widehat N_C(x_t)\)}
        \STATE \textbf{return} \(x_t\).
    \ELSIF[Smooth branch]{\(x_t\in\Omega\)}
        \STATE Set \(g_t=P_{T_C(x_t)}\nabla f(x_t)\), and choose
        \(s_t\in\mathrm{LMO}_{\|\cdot\|}(g_t)\).
        \STATE Choose the smallest \(j\in\mathbb N\) such that
        \(\mathcal R_{x_t}(\bar\alpha_M\beta_M^j s_t)\) is defined and
        \begin{equation}\label{eq:mcsd_pgd_smooth_branch}
        \begin{aligned}
            \alpha_t&=\bar\alpha_M\beta_M^j,
            &x_{t+1}&=\mathcal R_{x_t}(\alpha_t s_t),\\
            &&f(x_{t+1})&\le f(x_t)-c_M\alpha_t\|g_t\|_*.
        \end{aligned}
        \end{equation}
    \ELSE[PGD fallback]
        \STATE Choose the smallest \(j\in\mathbb N\) and a corresponding
        \(x_{t+1}\) satisfying
        \begin{equation}\label{eq:pgd_fallback}
        \begin{aligned}
            \alpha_t&=\bar\alpha_P\beta_P^j,
            &x_{t+1}
            &\in P_C\bigl(x_t-\alpha_t\nabla f(x_t)\bigr),\\
            &&f(x_{t+1})
            &\le f(x_t)
            +c_P\langle\nabla f(x_t),x_{t+1}-x_t\rangle.
        \end{aligned}
        \end{equation}
    \ENDIF
\ENDFOR
\end{algorithmic}
\end{algorithm}

\begin{remark}[Relation to P2GD--PGD]
For the Euclidean LMO, the smooth direction is
\(s_t=-g_t/\|g_t\|_{\mathbb E}\). With the metric-projection return map,
setting \(\tau_j:=\bar\alpha_M\beta_M^j/\|g_t\|_{\mathbb E}\) makes each
smooth trial and its Armijo test coincide with those of P2GD
\cite{SchneiderUschmajew2015P2GD,OlikierGallivanAbsil2023Stratified,olikier2026low}.
Thus, with the same branch and projection selections and this transformed
step-size grid, the Euclidean specialization of
Algorithm~\ref{alg:mcsd_pgd} agrees with P2GD--PGD up to normalization
\cite[Algorithms~5.1 and~7.2]{olikier2026low}. Non-Euclidean LMOs or other
return maps generally yield different smooth updates.
\end{remark}

\section{Convergence analysis}
\label{sec:convergence}

This section establishes convergence guarantees for the algorithms introduced
in Section~\ref{sec:MCSD}. For the \(C^2\) embedded-manifold case
\(C=\mathcal M\), we derive a best-iterate stationarity bound for
Algorithms~\ref{alg:MCSD} and~\ref{alg:MCSDTP}. For a general closed set
\(C\), we prove that Algorithm~\ref{alg:mcsd_pgd} is well defined and that
every accumulation point of its iterates, if any, is B-stationary.

\subsection{MCSD and MCSD-TP on embedded manifolds}
\label{subsec:mcsd_smooth_convergence}

Let \(\mathcal M\subset\mathbb E\) be a \(C^2\) embedded manifold. Under
either canonical return map in \eqref{eq:return_maps},
Algorithms~\ref{alg:MCSD} and~\ref{alg:MCSDTP} admit the common
representation
\[
    s_t\in\mathrm{LMO}_{\|\cdot\|}
    (\nabla_{\mathcal M}f(x_t)),
    \qquad
    x_{t+1}=\mathcal R_{x_t}(\alpha_t s_t).
\]
We therefore analyze both methods through the following quadratic return-map
model.
\begin{assumption}[Uniform return-map model]
\label{assumption:return_smooth}
Let \(f:\mathbb E\to\mathbb R\) be differentiable, and let
\(\mathcal R:\mathcal D\subset\mathcal M\times\mathbb E\to\mathcal M\)
be an ambient return map on \(\mathcal M\). There exist \(r,L>0\) such that,
for every \(x\in\mathcal M\) and \(u\in\mathbb E\) with \(\|u\|\le r\),
\((x,u)\in\mathcal D\) and
\begin{equation}\label{eq:return_map_descent}
f(\mathcal R_x(u))
    \le
    f(x)+\langle\nabla_{\mathcal M}f(x),u\rangle
    +\frac{L}{2}\|u\|_{\mathbb E}^2.
\end{equation}
\end{assumption}

The model radius is measured in the LMO norm, while the quadratic remainder
uses the Euclidean reference norm. Unless specified otherwise,
\(B(z,\rho):=\{y\in\mathbb E:\|y-z\|\le\rho\}\) denotes the closed ball in
the LMO norm. A uniform Lipschitz bound on the pullback gradients provides a
convenient sufficient condition for Assumption~\ref{assumption:return_smooth}.

\begin{lemma}[Uniformly Lipschitz pullback gradients]
\label{lemma:descent_lemma}
Let \(f:\mathbb E\to\mathbb R\) be differentiable, let
\(\mathcal R:\mathcal D\subset\mathcal M\times\mathbb E\to\mathcal M\)
be an ambient return map on \(\mathcal M\), and define
\(\phi_x:=f\circ\mathcal R_x\). Suppose there exist \(r,L>0\) such that
\(\mathcal M\times B(0,r)\subset\mathcal D\) and
\[
    \|\nabla_u\phi_x(u)-\nabla_u\phi_x(v)\|_{\mathbb E}
    \le L\|u-v\|_{\mathbb E}
\]
for every \(x\in\mathcal M\) and \(u,v\in B(0,r)\). Then
\(\mathcal R\) satisfies Assumption~\ref{assumption:return_smooth}.
\end{lemma}

\begin{proof}
Definition~\ref{def:ambient_return_map} and the chain rule give
\[
    \phi_x(0)=f(x),
    \qquad
    \nabla_u\phi_x(0)
    =D\mathcal R_x(0)^*\nabla f(x)
    =P_{T_x\mathcal M}^*\nabla f(x)
    =\nabla_{\mathcal M}f(x).
\]
Since \(B(0,r)\) is convex, a Taylor bound on \(\phi_x\)
yields \eqref{eq:return_map_descent} and proves the claim.
\end{proof}

On a compact \(C^2\) embedded manifold, an ambient return map and a locally
Lipschitz objective gradient yield this uniform pullback-gradient bound.

\begin{proposition}[Compact return-map criterion]
\label{prop:check_assumption}
Let \(f:\mathbb E\to\mathbb R\) be \(C^1\) with locally Lipschitz gradient
near a compact \(C^2\) embedded manifold \(\mathcal M\), and let
\(\mathcal R:\mathcal D\subset\mathcal M\times\mathbb E\to\mathcal M\)
be an ambient return map on \(\mathcal M\). Then
Assumption~\ref{assumption:return_smooth} holds.
\end{proposition}

\begin{proof}
Since \(\mathcal M\) is compact and \(\mathcal D\) is open, there exists
\(r>0\) such that
\(\mathcal K:=\mathcal M\times B(0,r)\subset\mathcal D\). For
\(x\in\mathcal M\) and \(u\in B(0,r)\), write
\(A_u:=D_u\mathcal R(x,u)\) and
\(g_u:=\nabla f(\mathcal R_x(u))\). After decreasing \(r\) if necessary,
compactness gives uniform bounds \(B_1,B_2,G\) for
\(\|A_u\|_{\rm op}\), \(\|D^2_{uu}\mathcal R(x,u)\|_{\rm op}\), and
\(\|g_u\|_{\mathbb E}\), respectively. The local Lipschitz continuity of
\(\nabla f\), again combined with compactness, gives a common constant
\(L_f\) such that
\(\|g_u-g_v\|_{\mathbb E}\le
L_f\|\mathcal R_x(u)-\mathcal R_x(v)\|_{\mathbb E}\).
The bounds on the first and second derivatives of \(\mathcal R\) also give
\(\|\mathcal R_x(u)-\mathcal R_x(v)\|_{\mathbb E}
\le B_1\|u-v\|_{\mathbb E}\) and
\(\|A_u-A_v\|_{\rm op}\le B_2\|u-v\|_{\mathbb E}\). Since the chain rule
gives \(\nabla_u\phi_x(u)=A_u^*g_u\), it follows that
\[
    \|\nabla_u\phi_x(u)-\nabla_u\phi_x(v)\|_{\mathbb E}
    \le (B_1^2L_f+B_2G)\|u-v\|_{\mathbb E}.
\]
Lemma~\ref{lemma:descent_lemma} now proves the claim.
\end{proof}

For a compact \(C^3\) embedded manifold, the projection return map
\(\mathcal R^{\rm P}\) is \(C^2\) on a uniform neighborhood
\cite[Lemma~4]{absil2012projection}, and a \(C^2\) local retraction induces
the \(C^2\) return map \(\mathcal R^{\rm R}\) in
\eqref{eq:return_maps}. Consequently, if \(\nabla f\) is locally Lipschitz
near \(\mathcal M\), Proposition~\ref{prop:check_assumption}
verifies Assumption~\ref{assumption:return_smooth} for both canonical return
maps. The following theorem gives the resulting common best-iterate
stationarity guarantee for Algorithms~\ref{alg:MCSD}
and~\ref{alg:MCSDTP}. Its proof follows the descent framework used
in recent work on Muon \cite{shen2025convergence}, and the resulting bound is
consistent with standard nonconvex first-order guarantees; see, e.g.,
\cite{bottou2018optimization,ghadimi2013stochastic}. Since \(\mathbb E\) is finite-dimensional, the norm-equivalence constant
\(N:=\sup_{\|z\|\le1}\|z\|_{\mathbb E}\) is finite and satisfies
\(\|z\|_{\mathbb E}\le N\|z\|\) for every \(z\in\mathbb E\).

\begin{theorem}[Best-iterate stationarity of MCSD and MCSD-TP]
\label{thm:convergence_deterministic}
Let \(\mathcal R\) satisfy
Assumption~\ref{assumption:return_smooth} with constants \(r,L>0\).
Starting from \(x_0\in\mathcal M\), generate \((x_t)_{t\ge0}\) by
\[
    s_t\in\mathrm{LMO}_{\|\cdot\|}(\nabla_{\mathcal M}f(x_t)),\qquad
    x_{t+1}:=\mathcal R_{x_t}(\alpha_ts_t),
    \qquad t=0,1,\ldots .
\]
Let \(T\ge1\). If \(0<\alpha_t\le r\) for
\(t=0,\ldots,T-1\), then
\[
    \sum_{t=0}^{T-1}\alpha_t
    \|\nabla_{\mathcal M}f(x_t)\|_*
    \le
    f(x_0)-f(x_T)
    +\frac{L N^2}{2}
    \sum_{t=0}^{T-1}\alpha_t^2 .
\]
Moreover, suppose that \(f\) is bounded below on \(\mathcal M\), and let
\(\Delta_f\ge0\) satisfy
\(\Delta_f\ge f(x_0)-\inf_{x\in\mathcal M}f(x)\). For
\(\alpha_t=\eta/\sqrt{t+1}\), \(t=0,\ldots,T-1\), where
\(0<\eta\le r\), we have
\[
    \min_{0\le t<T}\|\nabla_{\mathcal M}f(x_t)\|_*
    \le
    \frac{\displaystyle
        \Delta_f+\frac{L N^2\eta^2}{2}(1+\log T)}
    {\eta\sqrt T}
    =
    O\!\left(\frac{\log T}{\sqrt T}\right).
\]
\end{theorem}

\begin{proof}
Let $g_t:= \nabla_{\mathcal M}f(x_t)$. The LMO identity and the definition of \(N\) give
\[
    \langle g_t,s_t\rangle=-\|g_t\|_*,
    \qquad
    \|s_t\|_{\mathbb E}\le N .
\]
Since \(\|\alpha_ts_t\|\le\alpha_t\le r\),
Assumption~\ref{assumption:return_smooth} yields
\[
    f(x_{t+1})
    \le
    f(x_t)-\alpha_t\|g_t\|_*
    +\frac{L N^2}{2}\alpha_t^2 .
\]
Summing from \(t=0\) to \(T-1\) gives
\[
    \sum_{t=0}^{T-1}\alpha_t\|g_t\|_*
    \le
    f(x_0)-f(x_T)
    +\frac{L N^2}{2}
    \sum_{t=0}^{T-1}\alpha_t^2 .
\]
Now let \(\alpha_t=\eta/\sqrt{t+1}\). Since
\(f(x_T)\ge\inf_{\mathcal M}f\), the telescoping bound gives
\[
    \eta\sum_{t=0}^{T-1}\frac{\|g_t\|_*}{\sqrt{t+1}}
    \le
    \Delta_f+\frac{L N^2\eta^2}{2}
    \sum_{t=1}^T\frac1t .
\]
Therefore,
\[
    \min_{0\le t<T}\|g_t\|_*
    \le
    \frac{\displaystyle
        \Delta_f+\frac{L N^2\eta^2}{2}
        \sum_{t=1}^Tt^{-1}}
    {\displaystyle \eta\sum_{t=1}^Tt^{-1/2}}.
\]
Finally, the estimates
\(\sum_{t=1}^Tt^{-1}\le1+\log T\) and
\(\sum_{t=1}^Tt^{-1/2}\ge\sqrt T\) give the result.
\end{proof}

\begin{remark}[Best-iterate stationarity of Manifold Muon]
The proof of Theorem~\ref{thm:convergence_deterministic} also applies to an
exact tangent-constrained LMO step. Manifold Muon \cite{bernstein2025manifolds,kexuefm-11221,buchanan2025mmuonadmm,cesista2025spectralclipping,solonko2026skewon} uses the projection
update~\eqref{eq:proj_update} with \(g_t:=\nabla_{\mathcal M}f(x_t)\) and
the tangent-constrained LMO output
\begin{equation}\label{eq:manifold_muon}
    s_t\in\operatorname*{argmin}
    \{\langle g_t,s\rangle:\|s\|\le1,\ s\in T_{x_t}\mathcal M\}.
\end{equation}
If \(g_t\ne0\), then \(-g_t/\|g_t\|\) is feasible in
\eqref{eq:manifold_muon}; hence, by equivalence of norms in
\(\mathbb E\), there exists \(\kappa>0\), depending only on
\(\|\cdot\|\) and \(\|\cdot\|_{\mathbb E}\), such that
\[
    \langle g_t,s_t\rangle
    \le-\frac{\|g_t\|_{\mathbb E}^2}{\|g_t\|}
    \le-\kappa\|g_t\|_*.
\]
When \(g_t=0\), the final inequality
\(\langle g_t,s_t\rangle\le-\kappa\|g_t\|_*\) holds trivially. Moreover,
the constraint in \eqref{eq:manifold_muon} gives
\(\|s_t\|_{\mathbb E}\le N\). Thus, if the projection return map satisfies
Assumption~\ref{assumption:return_smooth} and \(0<\alpha_t\le r\), then
\eqref{eq:return_map_descent} gives
\[
    f(x_{t+1})
    \le f(x_t)-\kappa\alpha_t\|g_t\|_*
    +\frac{L N^2}{2}\alpha_t^2.
\]
The telescoping argument in Theorem~\ref{thm:convergence_deterministic}
therefore yields the same best-iterate bound, multiplied by
\(1/\kappa\).
\end{remark}

\subsection{MCSD--PGD on closed sets}
\label{subsec:mcsd_pgd_theory}
We now turn to the general closed-set problem
\eqref{eq:closed_problem_intro}, where \(C\subset\mathbb E\) is nonempty and
closed and \(f:\mathbb E\to\mathbb R\) is \(C^1\). Recall from
Section~\ref{sec:closed_set_prelim} that the set of B-stationary points is
\[
    S_B:=\{x\in C:-\nabla f(x)\in\widehat N_C(x)\}.
\]
We will prove that every iteration of MCSD--PGD is well defined
and that every accumulation point of its iterates, if any, is B-stationary.
We begin by relating B-stationarity to the intrinsic tangent-space gradients
on the safe region from Definition~\ref{def:safe_region}.

\begin{proposition}[Safe-region stationarity]
\label{prop:safe_region_stationarity}
Let \(C\subset\mathbb E\) be nonempty and closed, let
\(f:\mathbb E\to\mathbb R\) be \(C^1\), and let
\((\Omega,\mathcal R)\) be a safe region with return maps. For
\(x\in\Omega\), define
\[
    g(x):=P_{T_C(x)}\nabla f(x).
\]
Then \(T_C(x)\) is a linear subspace for every \(x\in\Omega\), so this
definition is well posed, and:
\begin{enumerate}
    \item[\textup{(i)}] if \(x\in\Omega\) and \(\mathcal M\) is any local
    manifold representative of \(C\) at \(x\), then
    \[
        T_C(x)=T_x\mathcal M,
        \qquad
        \widehat N_C(x)=T_C(x)^\perp;
    \]
    \item[\textup{(ii)}] for every \(x\in\Omega\),
    \(\|g(x)\|_*=0\) if and only if \(x\in S_B\);
    \item[\textup{(iii)}] if \(x_j\in\Omega\), \(x_j\to\bar x\in C\),
    and \(\|g(x_j)\|_*\to0\), then
    \(\bar x\in S_B\);
    \item[\textup{(iv)}] if \(\bar x\in C\setminus S_B\), then there exist
    a neighborhood \(U\) of \(\bar x\) and \(\eta>0\) such that
    \[
        \|g(x)\|_*\ge\eta
        \qquad\text{for all }x\in U\cap\Omega.
    \]
\end{enumerate}
\end{proposition}

\begin{proof}
\emph{Part \textup{(i)}.}
Fix \(x\in\Omega\), and choose an open neighborhood \(U_x\) and an embedded
\(C^3\) manifold \(\mathcal M_x\) such that
\(C\cap U_x=\mathcal M_x\cap U_x\). Since the tangent and regular normal
cones are local,
\[
    T_C(x)=T_x\mathcal M_x,
    \qquad
    \widehat N_C(x)
    =\widehat N_{\mathcal M_x}(x)
    =(T_x\mathcal M_x)^\perp.
\]
In particular, \(T_C(x)\) is a linear subspace. The same argument applies
to every local manifold representative of \(C\) at \(x\).

\emph{Part \textup{(ii)}.}
By part~\textup{(i)}, \(\|g(x)\|_*=0\) if and only if
\(\nabla f(x)\in T_C(x)^\perp=\widehat N_C(x)\). Since
\(\widehat N_C(x)\) is a linear subspace at every \(x\in\Omega\), this is
equivalent to \(-\nabla f(x)\in\widehat N_C(x)\), and hence to
\(x\in S_B\).

\emph{Part \textup{(iii)}.}
Because \(\Omega\) is relatively closed in \(C\), the assumptions imply
\(\bar x\in\Omega\). Choose an open neighborhood \(U_{\bar x}\) and an
embedded \(C^3\) manifold \(\mathcal M_{\bar x}\) such that
\[
    C\cap U_{\bar x}
    =\mathcal M_{\bar x}\cap U_{\bar x}.
\]
For all sufficiently large \(j\), one has \(x_j\in U_{\bar x}\), and
locality of the tangent cone gives
\[
    T_C(x_j)=T_{x_j}\mathcal M_{\bar x},
    \qquad
    T_C(\bar x)=T_{\bar x}\mathcal M_{\bar x}.
\]
The tangent-space projector varies continuously on
\(\mathcal M_{\bar x}\), and \(\nabla f\) is continuous. Consequently,
\(g(x_j)\to g(\bar x)\). Since \(\|g(x_j)\|_*\to0\), one has
\(g(\bar x)=0\), and part~\textup{(ii)} yields \(\bar x\in S_B\).

\emph{Part \textup{(iv)}.}
Suppose that the conclusion fails at some
\(\bar x\in C\setminus S_B\). Then, for every \(j\ge1\), there exists
\(x_j\in B(\bar x,1/j)\cap\Omega\) such that
\(\|g(x_j)\|_*<1/j\). Hence \(x_j\to\bar x\) and
\(\|g(x_j)\|_*\to0\). Part~\textup{(iii)} gives
\(\bar x\in S_B\), a contradiction.
\end{proof}

The preceding result uses only the safe-region geometry in
Definition~\ref{def:safe_region} and the \(C^1\) regularity of \(f\).
For the convergence analysis of Algorithm~\ref{alg:mcsd_pgd}, we additionally
assume that \(\nabla f\) is locally Lipschitz, which yields the following
local quadratic model.

\begin{lemma}[Local quadratic return-map model]
\label{lemma:local_quadratic_model}
Let \(C\subset\mathbb E\) be nonempty and closed, let
\(f:\mathbb E\to\mathbb R\) be \(C^1\) with locally Lipschitz gradient, and
let \((\Omega,\mathcal R)\) be a safe region with return maps in the sense of Definition~\ref{def:safe_region}. For
\(x\in\Omega\), write \(g(x):=P_{T_C(x)}\nabla f(x)\). For every
\(\bar x\in\Omega\), there exist a neighborhood \(U_{\bar x}\) of
\(\bar x\) and constants \(\tau_{\bar x},L_{\bar x}>0\) such that, for
every \(x\in U_{\bar x}\cap\Omega\) and \(u\in\mathbb E\) with
\(\|u\|\le\tau_{\bar x}\), one has \(u\in\mathcal D_x\) and
\[
    f(\mathcal R_x(u))
    \le f(x)+\langle g(x),u\rangle
       +\frac{L_{\bar x}}{2}\|u\|_{\mathbb E}^2.
\]
\end{lemma}

\begin{proof}
Fix \(\bar x\in\Omega\). Definition~\ref{def:safe_region} supplies a
neighborhood \(U\) on which \(C\) agrees with an embedded manifold
\(\mathcal M\) and the maps \(\mathcal R_x\) agree locally with the fibers
of an ambient return map
\(\widetilde{\mathcal R}:\widetilde{\mathcal D}\to\mathcal M\).
Choose a neighborhood \(U_{\bar x}\) of \(\bar x\) whose compact closure
lies in \(U\). Since \(C\) is closed and \(\Omega\) is relatively closed in
\(C\), the set \(K:=\Omega\cap\overline{U_{\bar x}}\) is compact.
Openness of \(\widetilde{\mathcal D}\), compactness of \(K\), and norm
equivalence give \(\tau_{\bar x}>0\) such that
\((x,u)\in\widetilde{\mathcal D}\) whenever \(x\in K\) and
\(\|u\|\le\tau_{\bar x}\). On this compact cylinder, the joint \(C^2\)
regularity of \(\widetilde{\mathcal R}\), local Lipschitz continuity of
\(\nabla f\), and the chain rule give \(L_{\bar x}>0\) such that, for
\(\phi_x:=f\circ\mathcal R_x\),
\[
    \|\nabla\phi_x(u)-\nabla\phi_x(v)\|_{\mathbb E}
    \le L_{\bar x}\|u-v\|_{\mathbb E}.
\]
Here \(x\) ranges over the smaller neighborhood in \(\Omega\), and
\(\|u\|,\|v\|\le\tau_{\bar x}\). Local equality with \(\mathcal M\) and
the ambient-return property give
\(\nabla\phi_x(0)=P_{T_C(x)}\nabla f(x)=g(x)\). The Taylor argument in
Lemma~\ref{lemma:descent_lemma} now yields the inequality in
the statement.
\end{proof}

An iterate outside the safe region is handled by the PGD
branch of Algorithm~\ref{alg:mcsd_pgd}. The next lemma shows that this branch
gives a uniform decrease near any non-B-stationary point, as a consequence
of the PGD analysis in \cite{OlikierWaldspurger2025PGD}.

\begin{lemma}[Uniform decrease of the PGD fallback]
\label{lemma:pgd_fallback_decrease}
Let \(C\subset\mathbb E\) be nonempty and closed, and let
\(f:\mathbb E\to\mathbb R\) be \(C^1\). Fix a safe region with return maps
\((\Omega,\mathcal R)\) and the PGD line-search parameters in
Algorithm~\ref{alg:mcsd_pgd}. For every
\(\bar x\in C\setminus S_B\), there exist a neighborhood \(U_P\) of
\(\bar x\) and \(\delta_P>0\) such that, for every iterate
\(x_t\in (U_P\cap C)\setminus(\Omega\cup S_B)\), the line
search~\eqref{eq:pgd_fallback} terminates after finitely many reductions
and its accepted iterate satisfies
\[
    f(x_{t+1})\le f(x_t)-\delta_P .
\]
\end{lemma}

\begin{proof}
Apply \cite[Proposition~5.3]{OlikierWaldspurger2025PGD} with
\(\beta=\beta_P\), \(c=c_P\), and initial step \(\bar\alpha_P\). That
proposition yields a neighborhood \(U_P\) of \(\bar x\) and
\(a_P\in(0,\bar\alpha_P]\) such that the Armijo condition in
\eqref{eq:pgd_fallback} holds for every
\(x_t\in(U_P\cap C)\setminus(\Omega\cup S_B)\), every projected point, and
every \(\alpha\in[a_P,a_P/\beta_P]\). The first point of the grid
\((\bar\alpha_P\beta_P^j)_{j\in\mathbb N}\) not exceeding
\(a_P/\beta_P\) lies in this interval. Hence the line search
\eqref{eq:pgd_fallback} terminates with
\(\alpha_t\in[a_P,\bar\alpha_P]\).

We next show that, after possibly shrinking \(U_P\), the accepted
displacements are uniformly bounded away from zero. Otherwise, there would
exist iterates
\(x_j\in C\setminus(\Omega\cup S_B)\) with \(x_j\to\bar x\), accepted step
sizes
\(\alpha_j\in[a_P,\bar\alpha_P]\), and accepted projection points \(x_j^+\)
such that
\[
    x_j^+\in P_C(x_j-\alpha_j\nabla f(x_j)),
    \qquad \|x_j^+-x_j\|_{\mathbb E}\to0.
\]
Along a subsequence, \(\alpha_j\to\alpha_*>0\). Closedness of the projection
graph then gives
\(\bar x\in P_C(\bar x-\alpha_*\nabla f(\bar x))\). By the projection
characterization of proximal normals
\cite[Example~6.16]{RockafellarWets1998}, this makes
\(-\nabla f(\bar x)\) a proximal, hence regular, normal to \(C\) at
\(\bar x\), contradicting \(\bar x\notin S_B\). Hence there exists
\(\rho_P>0\) such that
\(\|x_{t+1}-x_t\|_{\mathbb E}\ge\rho_P\) at every such iteration.

Finally, the defining projection inequality
\cite[Proposition~2.1]{OlikierWaldspurger2025PGD} gives
\begin{equation}\label{eq:pgd_projection_inequality}
    \langle\nabla f(x_t),x_{t+1}-x_t\rangle
    \le-\frac{\|x_{t+1}-x_t\|_{\mathbb E}^2}{2\alpha_t}
    \le-\frac{\rho_P^2}{2\bar\alpha_P}.
\end{equation}
The Armijo condition \eqref{eq:pgd_fallback} proves the result with
\(\delta_P=c_P\rho_P^2/(2\bar\alpha_P)\).
\end{proof}

We now combine Lemma~\ref{lemma:local_quadratic_model} with the PGD fallback
lemma to obtain the main closed-set guarantee. The proof parallels the
analysis of P2GD--PGD in \cite[Section~7]{olikier2026low}.

\begin{theorem}[B-stationary accumulation of MCSD--PGD]
\label{thm:mcsd_pgd_bstationary}
Let \(C\subset\mathbb E\) be nonempty and closed, and let
\(f:\mathbb E\to\mathbb R\) be \(C^1\) with locally Lipschitz gradient.
Let \((\Omega,\mathcal R)\) be a safe region with return maps in the sense of Definition~\ref{def:safe_region}.
Starting from any \(x_0\in C\), run Algorithm~\ref{alg:mcsd_pgd} with
parameters satisfying the requirements stated there.

At every nonterminal iteration of Algorithm~\ref{alg:mcsd_pgd}, the
applicable line search---\eqref{eq:mcsd_pgd_smooth_branch} on \(\Omega\) or
\eqref{eq:pgd_fallback} outside \(\Omega\)---terminates after finitely many
backtracking reductions and produces \(x_{t+1}\in C\).
Consequently, the algorithm either returns a point in \(S_B\) after
finitely many iterations or generates an infinite sequence
\((x_t)_{t\in\mathbb N}\), every accumulation point of which, if any,
belongs to \(S_B\).
\end{theorem}

\begin{proof}
\emph{Step 1: Finite termination and feasibility.}
Consider a nonterminal iteration, so \(x_t\notin S_B\). If
\(x_t\in\Omega\), then
Proposition~\ref{prop:safe_region_stationarity}\textup{(ii)} gives
\(\|g_t\|_*>0\). Lemma~\ref{lemma:local_quadratic_model}, applied at
\(\bar x=x_t\), supplies constants \(\tau_t,L_t>0\). Recall that
\(\|z\|_{\mathbb E}\le N\|z\|\) for every \(z\in\mathbb E\). Consider any
\(\alpha\) satisfying
\begin{equation}\label{eq:smooth_branch_acceptance_threshold}
    0<\alpha\le
    \min\left\{\tau_t,
    \frac{2(1-c_M)\|g_t\|_*}{L_tN^2}\right\}.
\end{equation}
The first bound in \eqref{eq:smooth_branch_acceptance_threshold}, together
with \(\|s_t\|\le1\), places \(\alpha s_t\) in the local model domain, so
the trial point is defined and belongs to \(C\). The local model, the LMO
identity
\(\langle g_t,s_t\rangle=-\|g_t\|_*\), and
\(\|s_t\|_{\mathbb E}\le N\) then give
\begin{equation}\label{eq:smooth_branch_acceptance_estimate}
    f(\mathcal R_{x_t}(\alpha s_t))
    \le f(x_t)-\alpha\|g_t\|_*
       +\frac{L_tN^2}{2}\alpha^2
    \le f(x_t)-c_M\alpha\|g_t\|_* .
\end{equation}
The second bound in \eqref{eq:smooth_branch_acceptance_threshold} gives the
last inequality in \eqref{eq:smooth_branch_acceptance_estimate}. Thus every
\(\alpha\) satisfying \eqref{eq:smooth_branch_acceptance_threshold} meets
the feasibility and sufficient-decrease requirements in
\eqref{eq:mcsd_pgd_smooth_branch}.
Since \(\bar\alpha_M\beta_M^j\to0\), the line search
\eqref{eq:mcsd_pgd_smooth_branch} terminates after finitely many reductions.
If \(x_t\notin\Omega\), then
\(x_t\in C\setminus(\Omega\cup S_B)\), and
Lemma~\ref{lemma:pgd_fallback_decrease}, applied with \(\bar x=x_t\), gives
finite termination of \eqref{eq:pgd_fallback}. Thus every nonterminal
iteration is well defined and preserves feasibility.
If the algorithm returns, its stopping test
\(-\nabla f(x_t)\in\widehat N_C(x_t)\) shows that \(x_t\in S_B\).

\emph{Step 2: Uniform decrease near a non-B-stationary accumulation point.}
Suppose the algorithm generates an infinite sequence
\((x_t)_{t\in\mathbb N}\), let \(x_*\) be an accumulation point, and assume
for contradiction that \(x_*\notin S_B\). The stopping test is never met,
so \(x_t\notin S_B\) for all \(t\). Applying
Lemma~\ref{lemma:pgd_fallback_decrease} at \(\bar x=x_*\) therefore gives a
neighborhood \(U_P\) of \(x_*\) and \(\delta_P>0\) such that
\(f(x_{t+1})\le f(x_t)-\delta_P\) for every iterate
\(x_t\in U_P\cap(C\setminus\Omega)\). This estimate controls fallback
steps in both cases below.

We show that there exist a neighborhood \(U\) of \(x_*\) and
\(\delta>0\) such that
\(f(x_{t+1})\le f(x_t)-\delta\) whenever \(x_t\in U\).

\emph{Case 1: \(x_*\notin\Omega\).}
Since \(\Omega\) is relatively closed in \(C\),
there exists a neighborhood \(U_0\) of \(x_*\) such that
\(U_0\cap C\subset C\setminus\Omega\). Set
\(U:=U_0\cap U_P\) and \(\delta:=\delta_P\). Thus only the PGD rule
\eqref{eq:pgd_fallback} is used at iterates in \(U\), and the preceding
PGD estimate gives the desired uniform decrease.

\emph{Case 2: \(x_*\in\Omega\).}
By
Proposition~\ref{prop:safe_region_stationarity}\textup{(iv)}, there exist a
neighborhood \(U_g\) of \(x_*\) and \(\eta>0\) such that
\[
    \|g(x)\|_*\ge\eta
    \qquad (x\in U_g\cap\Omega).
\]
Lemma~\ref{lemma:local_quadratic_model}, applied at \(\bar x=x_*\), gives
a neighborhood \(U_R\) and constants \(\tau_{x_*},L_{x_*}>0\). Set
\(U_M:=U_g\cap U_R\) and define
\begin{equation}\label{eq:uniform_smooth_acceptance_threshold}
    A_M
    :=
    \min\left\{
        \tau_{x_*},
        \frac{2(1-c_M)\eta}{L_{x_*}N^2}
    \right\}
    >0 .
\end{equation}
For every \(x_t\in U_M\cap\Omega\), one has
\(g_t=g(x_t)\) and \(\|g_t\|_*\ge\eta\). Consequently,
\eqref{eq:uniform_smooth_acceptance_threshold} and the calculation in
\eqref{eq:smooth_branch_acceptance_estimate}, now using the uniform constants
\(\tau_{x_*},L_{x_*}\), show that every \(0<\alpha\le A_M\) is acceptable
in \eqref{eq:mcsd_pgd_smooth_branch}. If
\(\bar\alpha_M\le A_M\), the initial trial is accepted. Otherwise, the first
grid point not exceeding \(A_M\) is greater than \(\beta_MA_M\) and is
acceptable. The line search accepts no later than this trial, so its step
satisfies
\[
    \alpha_t\ge
    \underline\alpha_M:=
    \min\{\bar\alpha_M,\beta_M A_M\}>0 ,
\]
and consequently
\[
    f(x_{t+1})
    \le
    f(x_t)-c_M\underline\alpha_M\eta .
\]
Set
\[
    U:=U_M\cap U_P,
    \qquad
    \delta:=\min\{c_M\underline\alpha_M\eta,\delta_P\}>0 .
\]
For \(x_t\in U\cap\Omega\), the smooth rule
\eqref{eq:mcsd_pgd_smooth_branch} gives a decrease of at least
\(c_M\underline\alpha_M\eta\). For
\(x_t\in U\cap(C\setminus\Omega)\),
Lemma~\ref{lemma:pgd_fallback_decrease} gives a decrease of at least
\(\delta_P\). Thus the desired uniform decrease also holds when
\(x_*\in\Omega\).

\emph{Step 3: Contradiction.}
Choose a subsequence \(x_{t_j}\to x_*\) with strictly increasing indices.
For all large \(j\),
\(x_{t_j}\in U\), so Step~2 gives
\[
    f(x_{t_j+1})\le f(x_{t_j})-\delta .
\]
The smooth acceptance rule \eqref{eq:mcsd_pgd_smooth_branch} decreases
\(f\). For every PGD projection, the defining projection inequality
\cite[Proposition~2.1]{OlikierWaldspurger2025PGD}---the first inequality in
\eqref{eq:pgd_projection_inequality}---makes the inner product in
\eqref{eq:pgd_fallback} nonpositive. Hence \(f(x_t)\) is nonincreasing, so
\[
    f(x_{t_{j+1}})
    \le
    f(x_{t_j+1})
    \le
    f(x_{t_j})-\delta.
\]
Letting \(j\to\infty\) and using continuity of \(f\) yields
\(f(x_*)\le f(x_*)-\delta\), a contradiction. Therefore \(x_*\in S_B\).
\end{proof}
\section{Safe regions of common closed-set constraints}
\label{sec:examples_safe_regions}
This section records safe regions and return
maps for three common constraints. We state only the ingredients required by
Definition~\ref{def:safe_region}. For objectives with locally Lipschitz
gradients, Theorem~\ref{thm:mcsd_pgd_bstationary} guarantees the convergence of MCSD--PGD.

\begin{example}[Sparse vectors]
\label{ex:sparse_safe_region}
For \(1\le r\le n\), consider
\[
    C:=\{x\in\mathbb R^n:\|x\|_0\le r\}.
\]
Such cardinality constraints arise in compressed sensing
\cite{blumensath2009iterative,foucart2011hard}, best-subset regression
\cite{bertsimas2016best}, sparse principal component analysis
\cite{bertsimas2020solving}, and sparse phase retrieval
\cite{cai2016optimal,li2013sparse,
moravec2007compressive,shechtman2014gespar}. For \(\Delta>0\), define
\[
    \Omega_\Delta
    :=\left\{x\in C:\|x\|_0=r,\
      \min_{i\in\operatorname{supp}(x)}|x_i|\ge\Delta\right\}.
\]
For \(S\subset\{1,\ldots,n\}\), let
\(\mathbb R^S:=\{z\in\mathbb R^n:z_{S^c}=0\}\). Both \(C\) and
\(\Omega_\Delta\) are finite unions of closed subsets of these coordinate
subspaces, so \(\Omega_\Delta\) is relatively closed in \(C\).

For \(x\in\Omega_\Delta\), set \(S_x:=\operatorname{supp}(x)\) and define
\(\mathcal R_x(u):=x+P_{S_x}u\), \(u\in\mathbb R^n\), where \(P_{S_x}\)
is the Euclidean projector onto \(\mathbb R^{S_x}\). Fix
\(\bar x\in\Omega_\Delta\), let \(S:=S_{\bar x}\), and set
\(U_{\bar x}:=\{x:\|x-\bar x\|_{\mathbb E}<\Delta/2\}\). Then every
\(x\in\Omega_\Delta\cap U_{\bar x}\) has support \(S\), and
\[
    C\cap U_{\bar x}=\mathbb R^S\cap U_{\bar x}.
\]
The space \(\mathbb R^S\) is a smooth embedded manifold, and the map
\[
    \widetilde{\mathcal R}:
    (\mathbb R^S\cap U_{\bar x})\times\mathbb R^n\to\mathbb R^S,
    \qquad
    \widetilde{\mathcal R}(x,u):=x+P_Su,
\]
is jointly \(C^\infty\), with
\(\widetilde{\mathcal R}_x(0)=x\) and
\(D\widetilde{\mathcal R}_x(0)=P_S\). Thus
\((\Omega_\Delta,\mathcal R)\) satisfies
Definition~\ref{def:safe_region}.
\end{example}
In the next two matrix-constrained examples, let \(X\in\mathbb{R}^{m\times n}\), and let \(\sigma_k(X)\) denote its \(k\)th largest singular value for \(1\leq k\leq \min\{m,n\}\).
\begin{example}[Determinantal variety]
\label{ex:rank_constrained_safe_region}
For \(1\le r\le\min\{m,n\}\), consider the determinantal variety
\[
    C:=\{X\in\mathbb R^{m\times n}:\operatorname{rank}(X)\le r\}.
\]
Low-rank models of this form arise in matrix completion
\cite{vandereycken2013low}, recommender systems \cite{koren2009}, and
low-rank quantum-state tomography \cite{gross2010quantum}. For \(\Delta>0\), define
\[
    \mathcal M_r:=\{X\in\mathbb R^{m\times n}:\operatorname{rank}(X)=r\},
    \qquad
    \Omega_\Delta:=\{X\in C:\sigma_r(X)\ge\Delta\}.
\]
Continuity of \(\sigma_r\) makes \(\Omega_\Delta\) relatively closed in
\(C\), while \(\mathcal M_r\) is a \(C^\infty\) embedded manifold
\cite[Section~2]{vandereycken2013low}. For
\(\bar X\in\Omega_\Delta\), the Eckart--Young theorem
\cite{EckartYoung1936} gives
\(\operatorname{dist}_{\mathbb E}
(\bar X,\{Y:\operatorname{rank}(Y)\le r-1\})
=\sigma_r(\bar X)\ge\Delta\). Thus, for
\(U_{\bar X}:=\{Y:\|Y-\bar X\|_{\mathbb E}<\Delta/2\}\), $C\cap U_{\bar X}=\mathcal M_r\cap U_{\bar X}$.

Let \(\sigma_{r+1}:=0\) when \(r=\min\{m,n\}\), and define
\(\mathcal G_r:=\{Y:\sigma_r(Y)>\sigma_{r+1}(Y)\}\). The rank-\(r\)
truncated-SVD map \(\Pi_r\) is single-valued and \(C^\infty\) on
\(\mathcal G_r\). Define
\begin{equation}\label{eq:determinantal_return_maps}
    \begin{aligned}
    \mathcal D^{\rm P}
    &:=
    \{(X,u)\in\mathcal M_r\times\mathbb E:X+u\in\mathcal G_r\},
    &
    \mathcal R^{\rm P}(X,u)
    &:=
    \Pi_r(X+u),\\
    \mathcal D^{\rm R}
    &:=
    \{(X,u)\in\mathcal M_r\times\mathbb E:
      X+P_{T_X\mathcal M_r}u\in\mathcal G_r\},
    &
    \mathcal R^{\rm R}(X,u)
    &:=
    \Pi_r\bigl(X+P_{T_X\mathcal M_r}u\bigr).
\end{aligned}
\end{equation}
For a compact SVD \(X=U\Sigma V^\top\), the tangent projector acts on
\(H\in\mathbb E\) as
\(P_{T_X\mathcal M_r}H=UU^\top H+HVV^\top-UU^\top HVV^\top\).
Consequently, \(X+P_{T_X\mathcal M_r}H\) has rank at most \(2r\). Hence,
\(\mathcal R_X^{\rm R}(H)\) can be evaluated by truncating an SVD of a core of
size at most \(2r\times2r\), as in the P2GD branch of P2GD--PGD
\cite[Algorithms~5.1 and~7.2 and Section~8.1]{olikier2026low}. 
Both domains are open and contain \(\mathcal M_r\times\{0\}\), and both
maps are jointly \(C^\infty\) with values in \(\mathcal M_r\). The
derivative of the local metric projection
\cite[Lemma~4]{absil2012projection} gives
\[
    \mathcal R_X^{\rm P}(0)=\mathcal R_X^{\rm R}(0)=X,
    \qquad
    D\mathcal R_X^{\rm P}(0)
    =D\mathcal R_X^{\rm R}(0)
    =P_{T_X\mathcal M_r}.
\]
Thus either family equips \(\Omega_\Delta\) as a safe region with return
maps.
\end{example}

\begin{example}[Spectral-norm sphere]
\label{ex:spectral_sphere_safe_region}
Let \(R>0\) and \(\min\{m,n\}\ge2\), and consider the spectral-norm sphere
\[
    C:=\{X\in\mathbb R^{m\times n}:\|X\|_{2\to 2}=R\}.
\]
This set is nonempty and closed.
The spectral norm controls operator gain in Lipschitz networks
\cite{anil2019sorting}, while fixing it at a prescribed radius implements the
spectral scaling motivated by \(\mu\)P
\cite{yang2021tensorprogramsv,yang2023spectral,xie2026controlled}.

Fix \(\Delta\in(0,R)\), and define
\[
\begin{aligned}
    \mathcal V&:=\{Y:\sigma_1(Y)>\sigma_2(Y)\},
    &\mathcal M&:=C\cap\mathcal V,\\
    \Omega_\Delta&:=\{X\in C:\sigma_1(X)-\sigma_2(X)\ge\Delta\}.
\end{aligned}
\]
Continuity of the singular values makes \(\mathcal V\) open and
\(\Omega_\Delta\) relatively closed in \(C\). On \(\mathcal V\), the
function \(\sigma_1\) is \(C^\infty\), and
\(Q(Y):=\nabla\sigma_1(Y)\) satisfies \(\|Q(Y)\|_{\mathbb E}=1\)
\cite{Lewis1996Spectral,LewisSendov2005}. Thus \(\mathcal M\) is a
\(C^\infty\) embedded hypersurface and \(C\cap\mathcal V=\mathcal M\), so
\(\mathcal V\) supplies the required neighborhood of every point in
\(\Omega_\Delta\). Its tangent projector is
\[
    P_{T_X\mathcal M}H
    =H-\langle Q(X),H\rangle Q(X).
\]

For \(X\in\mathcal M\), set
\(Z_X(u):=X+P_{T_X\mathcal M}u\), and define
\[
    \mathcal D:=\{(X,u)\in\mathcal M\times\mathbb E:Z_X(u)\in\mathcal V\},
    \qquad
    \mathcal R(X,u):=R\,\frac{Z_X(u)}{\sigma_1(Z_X(u))}.
\]
The set \(\mathcal D\) is open in \(\mathcal M\times\mathbb E\) and
contains \(\mathcal M\times\{0\}\). Normalization sets the leading singular
value to \(R\) without changing the strict gap, so \(\mathcal R\) is jointly
\(C^\infty\) with values in \(\mathcal M\). Direct differentiation gives
\[
    \mathcal R_X(0)=X,
    \qquad
    D\mathcal R_X(0)=P_{T_X\mathcal M}.
\]
Restricting the fibers of \(\mathcal R\) to \(X\in\Omega_\Delta\) therefore
verifies Definition~\ref{def:safe_region}.
\end{example}

\section{Experiments}\label{sec:experiments}
We use three problems to assess the proposed methods and the influence of the
chosen LMO geometry: Stiefel-constrained PCA for MCSD, and weighted low-rank
approximation and sparse phase retrieval for MCSD--PGD. All experiments
are run on a MacBook Pro with an Apple M1 Max processor.

\subsection{Stiefel-constrained principal component analysis}
\label{sec:PCA}
\paragraph{Problem}
Given a data matrix \(X\in\mathbb R^{n\times d}\) and its Gram matrix
\(S:=XX^\top\), we consider the weighted PCA problem
\cite{wold1987principal}
\begin{equation}
    \min_{W\in\mathrm{St}(n,p)}
    f(W):=-\tfrac12\operatorname{tr}(W^\top SWD),
    \label{eq:weighted_pca}
\end{equation}
where \(D\in\mathbb R^{p\times p}\) is diagonal and nonnegative. This
gives the Brockett cost, a standard Stiefel benchmark
\cite{brockett1991dynamical,edelman1998geometry,absil2008optimization}.

\paragraph{Methods}
We compare MCSD under the five LMO norms
\(\|\cdot\|_{2\to 2}\), \(\|\cdot\|_F\),
\(\|\cdot\|_{1\to\infty}\), \(\|\cdot\|_{1\to 2}\), and
\(\|\cdot\|_{2\to\infty}\). The spectral-norm specialization
\(\|\cdot\|_{2\to 2}\) is SPEL, while the Frobenius-norm specialization
\(\|\cdot\|_F\) is NormRGD (Table~\ref{tab:mcsd-2x2}). We also compare SPEL
with four Manifold Muon variants that target the tangent-space spectral-norm
LMO: MM\textsubscript{Su}, MM\textsubscript{Bern}, and
MM\textsubscript{Buch} use, respectively, the fixed-point, dual-subgradient,
and ADMM solvers of
\cite{kexuefm-11221,bernstein2025manifolds,buchanan2025mmuonadmm}, whereas
MM\textsubscript{Skewon} uses the exact closed-form solver of
\cite{solonko2026skewon}.
\begin{figure}[!ht]
    \centering
    \includegraphics[width=\linewidth]{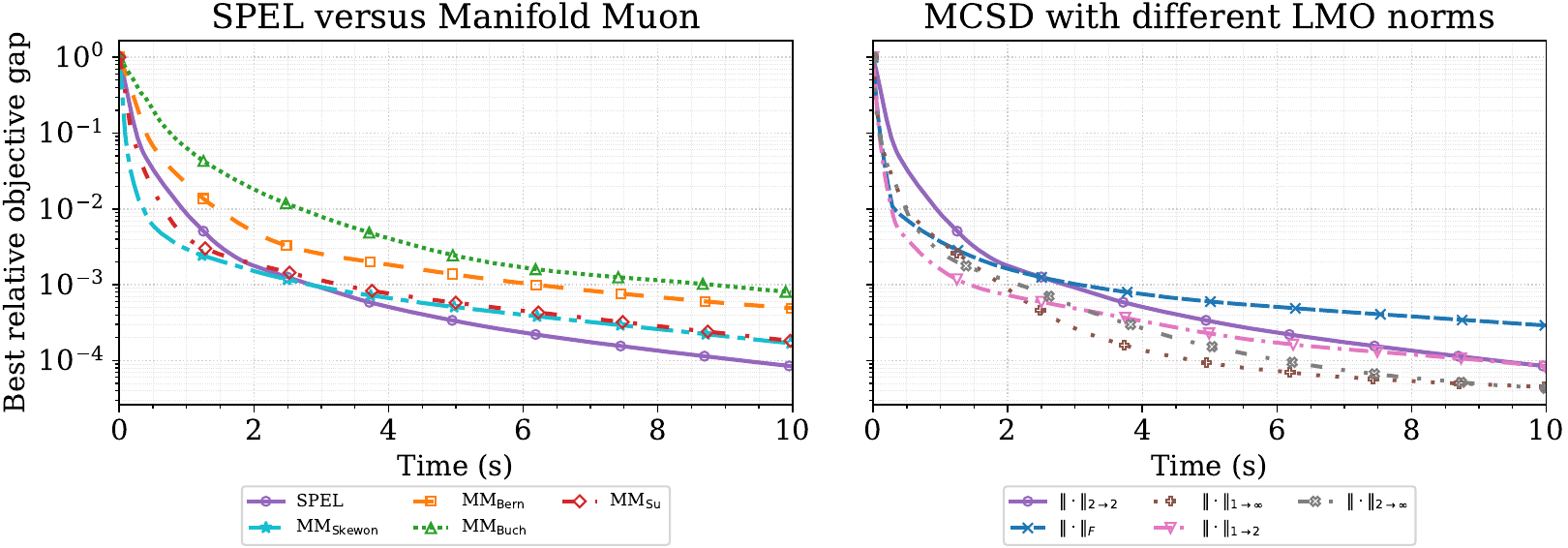}
    \caption{Best relative objective gap versus optimizer time for
    \((n,d,p)=(512,2048,32)\). The left panel compares SPEL with three
    iterative MM variants and exact MM\textsubscript{Skewon}; the right panel compares five LMO norms within MCSD.}
    \label{fig:pca_dis}
\end{figure}
\paragraph{Setup}
We use \((n,d,p)=(512,2048,32)\). To form a spiked instance, we scale the rows
of a standard Gaussian matrix so that the first \(p\) population eigenvalues
decrease geometrically from \(16\) to \(4\), while the remaining eigenvalues
equal \(1\); the diagonal entries of \(D\) decrease linearly from \(1\) to
\(0.2\). All methods use the same data, projected Gaussian initialization,
exact polar factors computed by SVD, and steps
\(\alpha_t=\eta/\sqrt{t+1}\). We tune each method independently over the same
26-point grid in \([10^{-3},10]\). Each candidate is evaluated in three
10-second runs, and the value with the smallest median final objective gap is
selected. The three iterative MM baselines use five inner iterations for
MM\textsubscript{Bern} and
MM\textsubscript{Buch} and two for MM\textsubscript{Su};
MM\textsubscript{Skewon} uses its exact compact solver. With the selected
\(\eta\) fixed, we perform five additional runs and report the trajectory
with median final objective gap. Figure~\ref{fig:pca_dis} plots
\(\min_{k:t_k\leq t}(f(W_k)-f^\star)/(f(W_0)-f^\star)\), where \(t_k\) is the
optimizer time after \(k\) updates and \(f^\star\) is obtained from the
eigendecomposition of \(S\).

\paragraph{Results}
Under the common tuning protocol, SPEL attains the smallest end-of-budget gap
among the Manifold Muon comparisons. The corresponding gaps are about twice
as large for MM\textsubscript{Skewon} and MM\textsubscript{Su}, six times as
large for MM\textsubscript{Bern}, and nine times as large for
MM\textsubscript{Buch}. The right panel shows that the LMO geometry materially
affects performance: the four non-Euclidean norms finish within a factor of
two and outperform the Frobenius-norm LMO by factors of about three to seven
on this instance.

\subsection{Weighted low-rank approximation}
\label{sec:low-rank-constrained-least-squares}
\paragraph{Problem}
We consider weighted low-rank approximation (WLRA) over
\(C=\{X\in\mathbb R^{m\times n}:\operatorname{rank}(X)\le r\}\):
\[
    \min_{X\in C} f(X):=\tfrac12\|\sqrt W\odot(X-A)\|_F^2,
\]
where \(W\) has positive entries, \(\odot\) denotes entrywise multiplication,
and \(A\in C\). Thus the minimum value is \(f^\star=0\). We use the apocalyptic construction of
\cite[Sections~8.2--8.3]{olikier2026low}, in which smooth rank-\(r\)
iterations may approach a nonstationary point of lower rank.

\paragraph{Methods}
We test two MCSD--PGD variants, distinguished by the LMO norm: the spectral
norm \(\|\cdot\|_{2\to2}\) and the weighted Frobenius norm
\(\|\cdot\|_{W,F}\). The latter is induced, as in
\eqref{eq:preconditioned_lmo}, by \(\mathcal B_W[D]:=W\odot D\). Since
\(\mathcal B_W[D]=\nabla^2 f(X)[D]\) is the constant ambient Hessian and
\(\mathcal B_W^{-1}\) is applied entrywise, this LMO incorporates exact
ambient inverse-curvature scaling at low cost.

Both variants use MCSD-TP on
\(\Omega_\Delta=\{X\in C:\sigma_r(X)\ge\Delta\}\), where \(C\) locally
coincides with the smooth rank-\(r\) manifold \(\mathcal M_r\). They project
the LMO output \(s\) onto \(T_X\mathcal M_r\) and return through the exact
small-core rank-\(r\) projection
\(\Pi_r(X+\alpha P_{T_X\mathcal M_r}s)\); see
Example~\ref{ex:rank_constrained_safe_region}. Outside \(\Omega_\Delta\), they
use PGD on \(C\). The baselines are PGD, P2GD--PGD
\cite[Algorithm~7.2]{olikier2026low}, RFDR \cite{OlikierAbsil2023RFDR}, and
the CRFDR-entry specialization of ERFDR \cite{OlikierAbsil2024ERFDR}; we use
the authors' implementations.\footnote{
\url{https://github.com/golikier/LowRankOptimization}}

\paragraph{Setup}
We use \((m,n,r,r')=(600,400,15,10)\), draw the entries of \(W\) uniformly
from \((0,1)\), and otherwise follow the construction and initialization in
\cite[Sections~8.2--8.3]{olikier2026low}. Every method is tuned independently
on instances disjoint from the \(100\) evaluation instances, using multistage
grid searches over its line-search parameters. 

For each MCSD--PGD variant, at a smooth-branch iteration let
\(d_t=P_{T_{X_t}\mathcal M_r}s_t\) and set
\(\eta_t:=\|g_t\|_*\). The smooth line
search starts from the
smallest of \(\bar\alpha_M\eta_t\), \(\bar\alpha_M\), and
\(\alpha_{M,t}^{\rm prev}/\beta_M\), and then backtracks by \(\beta_M\).
Here, \(\alpha_{M,t}^{\rm prev}\) is the step size accepted at the most recent
successful smooth-branch iteration.
Because P2GD and PGD use unnormalized gradient directions whereas an LMO
output is normalized in its prescribed geometry, \(\eta_t\) restores the
gradient scale in the initial trial. The previous accepted step provides a
local scale estimate while allowing the next trial to grow by one
backtracking level.
The common PGD fallback uses
\((\bar\alpha_P,\beta_P,c_P)=(3.2,1/3,0.0747664)\). The selected
\((\bar\alpha_M,\beta_M,c_M,\Delta)\) values are
\((4.8,0.5,3\times10^{-3},3\times10^{-3})\) for
\(\|\cdot\|_{2\to2}\) and
\((4.8,0.2,3\times10^{-4},3\times10^{-3})\) for
\(\|\cdot\|_{W,F}\). The P2GD--PGD tuple is
\((12.8,0.104,0.0379170,3\times10^{-4})\). The RFDR and CRFDR-entry tuples
are, respectively, \((51.2,0.0303030,0.1,0.32)\) and
\((51.2,0.0303030,0.0270270,0.0025)\). All methods use the same instances
and initializations and are run three times per instance, with a
\(60\)-second limit per run. 

\begin{table}[!ht]
\centering
\caption{Median WLRA hitting times to $f(X_t)\le \varepsilon$ in seconds over \(100\) random
instances. Each entry is the median of the per-instance median times among
instances solved in all three runs; solved counts are in parentheses. 
}
\label{tab:wlra-accuracy-speedup}
\setlength{\tabcolsep}{3.5pt}
\small
\begin{tabular}{@{}lrrrrrr@{}}
\toprule
& & & & & \multicolumn{2}{c}{MCSD--PGD} \\
\cmidrule(lr){6-7}
\(\varepsilon\) & PGD & P2GD--PGD & RFDR & CRFDR-entry & \(\|\cdot\|_{2\to2}\) & \(\|\cdot\|_{W,F}\) \\
\midrule
$10^{-3}$  & $.439\,(100)$ & $.0744\,(100)$ & $.250\,(100)$ & $.0752\,(100)$ & $.159\,(100)$ & $.0351\,(100)$ \\
$10^{-4}$  & $1.06\,(100)$ & $.137\,(100)$  & $.381\,(100)$ & $.115\,(100)$  & $.459\,(100)$ & $.0379\,(100)$ \\
$10^{-6}$  & $2.96\,(99)$  & $.308\,(100)$  & $.708\,(99)$  & $.245\,(100)$  & $1.97\,(99)$  & $.0388\,(100)$ \\
$10^{-8}$  & $4.35\,(97)$  & $.474\,(99)$   & $1.18\,(99)$  & $.440\,(99)$   & $3.85\,(97)$  & $.0418\,(100)$ \\
$10^{-10}$ & $5.81\,(96)$  & $.625\,(99)$   & $1.50\,(99)$  & $.675\,(99)$   & $6.33\,(95)$  & $.0448\,(100)$ \\
\bottomrule
\end{tabular}
\end{table}

\begin{figure}[!ht]
    \centering
    \includegraphics[width=\textwidth]{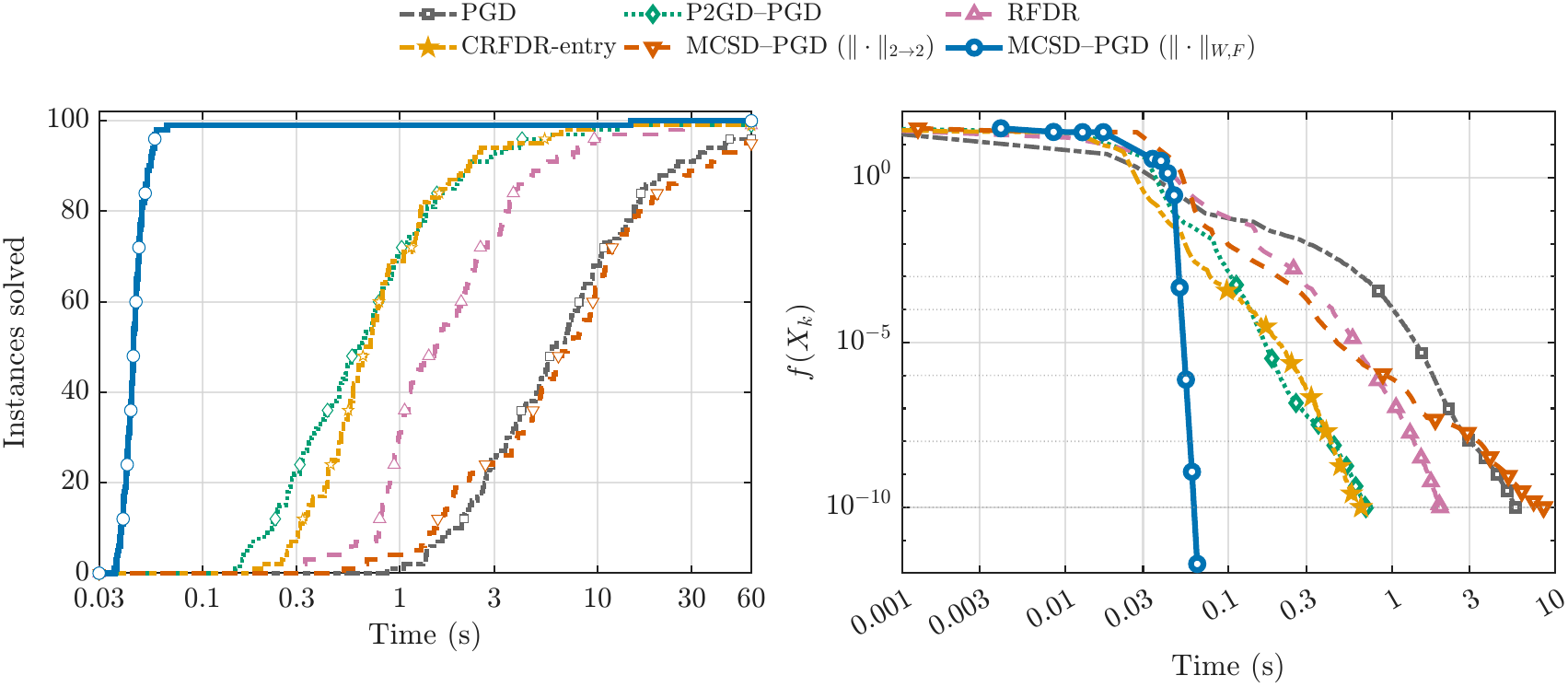}
    \caption{WLRA performance over \(100\) evaluation instances.
    Left: completion profile of the per-instance median time to reach
    \(10^{-10}\) across three runs. Right: objective histories plotted
    against median time across three runs on a
    representative instance.}
    \label{fig:wlra-performance}
\end{figure}

\paragraph{Results}
The results are summarized in Table \ref{tab:wlra-accuracy-speedup} and Figure \ref{fig:wlra-performance}. MCSD--PGD with \(\|\cdot\|_{W,F}\) has the smallest conditional median at every
tolerance and is the only method that solves all \(100\) instances at every
tolerance. It is \(3.61\) times faster than P2GD--PGD at
\(10^{-4}\) and \(13.96\) times faster at \(10^{-10}\). At the strictest
tolerance, P2GD--PGD, RFDR, and CRFDR-entry each solve \(99\)
instances, whereas PGD and spectral MCSD--PGD solve \(96\) and \(95\),
respectively. Spectral MCSD--PGD improves on PGD through \(10^{-8}\), but
not at \(10^{-10}\). These results suggest that the ambient preconditioning
enabled by MCSD can improve performance when the LMO geometry reflects the
problem curvature; they do not imply that arbitrary non-Euclidean LMOs
outperform gradient-based methods.

\subsection{Sparse phase retrieval}
\label{sec:sparse-phase-retrieval-experiment}
\paragraph{Problem}
We consider (real) sparse phase retrieval
\cite{moravec2007compressive,li2013sparse,shechtman2014gespar,cai2016optimal}:
\[
    \min_{\|x\|_0\le s} f(x):=\frac{1}{4m}\sum_{j=1}^m
    \big((a_j^\top x)^2-b_j\big)^2,
\]
where the clean observations satisfy \(b_j=(a_j^\top x^\natural)^2\) with $x^\natural\in \mathbb R^n$ and $\|x^\natural\|_0 \le s$. Thus the optimal value is \(f^\star=0\).

\paragraph{Methods}
We compare MCSD--PGD under three LMO norms: \(\|\cdot\|_2\),
\(\|\cdot\|_3\), and \(\|\cdot\|_\infty\). On the safe region
\(\Omega_\Delta\), each variant uses MCSD-TP with the sparse-vector return map
from Example~\ref{ex:sparse_safe_region}. Writing
\(S_t=\operatorname{supp}(x_t)\), we select an LMO output
\(s_t\in\mathbb R^{S_t}\), so the update reduces to
\(\mathcal R_{x_t}(\alpha_t s_t)=x_t+\alpha_tP_{S_t}s_t
=x_t+\alpha_t s_t\). Outside \(\Omega_\Delta\), the method falls back to a
hard-thresholded PGD step. The \(\|\cdot\|_2\) variant initializes
backtracking from a fixed trial step, whereas the other two use the adaptive
initialization from the WLRA experiment. The baselines are PGD \cite{blumensath2009iterative,foucart2011hard},
GraHTP--GN \cite{DaiLuYou2024GraHTP}, and Newton projection
\cite{CaiLongWenYing2024SparsePR}. For the latter two, we implement the
real-valued forms of the published updates; Newton projection follows the
paper's preconditioned conjugate-gradient convention.

\paragraph{Setup}
We use \((n,m,s)=(12000,8000,20)\), standard Gaussian sensing vectors, and
a unit-norm Gaussian signal on a uniformly sampled support. Every method starts
from the same support-blind sparse spectral initializer of
\cite[Algorithm~1]{JagatapHegde2019SampleEfficient}, also used in
\cite[Section~3.2]{CaiLongWenYing2024SparsePR}: it selects the \(s\) largest scores
\(m^{-1}\sum_j b_j a_{ji}^2\) and uses the scaled leading eigenvector of
\(m^{-1}\sum_j b_j a_{j,\widehat S_0}a_{j,\widehat S_0}^{\top}\) on the
selected support \(\widehat S_0\).

We scale the safe-region margin as
\(\Delta=\rho_\Delta\|x_0\|_2/\sqrt{s}\), so \(\rho_\Delta\) is a
dimensionless threshold controlling how far the active entries must remain
from zero for the smooth branch to be used.
For each LMO norm, we report \((\bar\alpha_M,\beta_M,\rho_\Delta)\), where
\(\bar\alpha_M\) caps the initial smooth step-size trial and \(\beta_M\) contracts a
rejected trial; see Algorithm \ref{alg:mcsd_pgd}. The selected values are \((0.6,0.25,0.02)\) for
\(\|\cdot\|_2\), \((3.2,0.25,0.005)\) for \(\|\cdot\|_3\), and
\((3.2,0.5,0.005)\) for \(\|\cdot\|_\infty\); all three use the smooth
Armijo constant \(c_M=10^{-3}\). The common PGD fallback uses
\((\bar\alpha_P,\beta_P,c_P)=(0.4,0.5,0.01)\), denoting its initial trial
step, backtracking factor, and Armijo constant, respectively. GraHTP--GN uses
support step \(0.4\) and one Gauss--Newton refit; Newton
projection uses \((\eta,\tau)=(0.75,1)\) and a conjugate-gradient tolerance
of \(10^{-3}\), with at most \(20\) inner iterations. We evaluate \(100\)
random instances and run each method three times under the same one-second
method-time budget. We use the best relative objective
\(\min_{k:t_k\le t} f(x_k)/f(x_0)\) as the optimality criterion. A run is solved when this quantity reaches \(10^{-10}\). An instance
enters the completion profile only if all three repeats solve, and its time
is the median of their hitting times.

\begin{figure}[!ht]
    \centering
    \includegraphics[width=\textwidth]{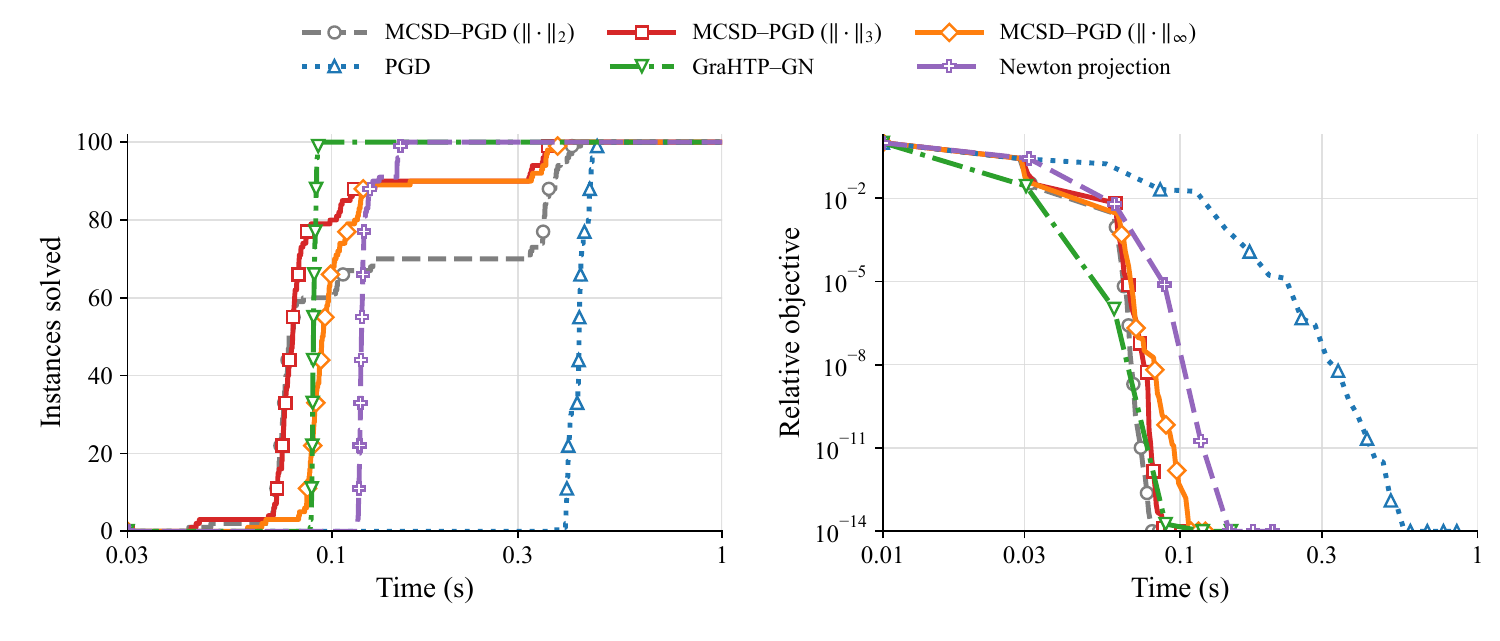}
    \caption{Sparse phase retrieval over \(100\) instances. Left: completion
    profile of the three-run median time to relative objective \(10^{-10}\);
    all runs must reach the cutoff. Right: a representative median trajectory.}
    \label{fig:spr-clean-performance}
\end{figure}

\paragraph{Results}
All methods solve all \(100\) instances (Figure~\ref{fig:spr-clean-performance}).
The median hitting times are \(0.0778\) s for
\(\|\cdot\|_2\) MCSD--PGD, \(0.0794\) s for \(\|\cdot\|_3\) MCSD--PGD,
\(0.0898\) s for GraHTP--GN, \(0.0948\) s for
\(\|\cdot\|_\infty\) MCSD--PGD, \(0.1192\) s for Newton projection, and
\(0.4299\) s for PGD. Thus all three MCSD--PGD variants are competitive with
the two specialized second-order methods and are between \(4.5\) and \(5.5\)
times faster than PGD.

\section{Concluding remarks}\label{sec:conclusion}
We introduced MCSD and MCSD-TP, two single-loop methods that combine
norm-induced ambient LMOs with projection or tangent projection and
retraction, thereby avoiding iterative tangent-space LMO subproblems. On
smooth embedded manifolds, both methods attain an
\(O(\log T/\sqrt{T})\) best-iterate stationarity bound. On closed sets with
certified smooth regions, MCSD--PGD combines these smooth steps with a PGD
fallback and guarantees that every accumulation point, if any, is Bouligand
stationary. The experiments on Stiefel-constrained PCA, weighted low-rank
approximation, and sparse phase retrieval indicate that the preferred search
geometry depends on the application; none of the tested norms or algorithms
is uniformly superior.

Several directions remain open. First, the current theory treats a fixed LMO
norm. Studying adaptive LMO geometries
\cite{becigneul2018riemannian,kasai2019riemannian} and comparing ambient with
Riemannian preconditioning \cite{mishra2016riemannian,bian2024preconditioned}
would clarify their computational tradeoffs.
Second, the present closed-set guarantee, like related accumulation-point
results \cite{OlikierGallivanAbsil2023Stratified,OlikierWaldspurger2025PGD}, is
asymptotic; an important question is whether uniform regularity or error-bound
assumptions yield a computable approximate Bouligand-stationarity measure and iteration complexity. Third, stochastic extensions of
MCSD--PGD, analogous to stochastic Riemannian methods
\cite{bonnabel2013stochastic}, could scale to applications such as
large-language-model (LLM) training on the spectral-norm sphere
\cite{yang2023spectral,xie2026controlled}
(Example~\ref{ex:spectral_sphere_safe_region}).

\bibliographystyle{alpha}
\bibliography{mybib}

\end{document}